\documentclass[11pt]{amsart}

\usepackage{amsmath}
\usepackage{amsfonts}
\usepackage{amssymb}
\usepackage{amsthm}
\usepackage{mathtools} % An improvement of amsmath
\usepackage{latexsym}

\usepackage{enumerate}
\usepackage{cancel}
\usepackage{cases}
\usepackage{empheq}
\usepackage{multicol}

\usepackage{subfigure}
\usepackage{float}

\usepackage{wrapfig}

\usepackage{color}

\usepackage[color=white,linecolor=black]{todonotes}
\usepackage{mathrsfs}
\usepackage{fontenc} %T1 font encoding
\usepackage{inputenc} %UTF-8 support

\usepackage{verbatim}

\theoremstyle{plain}
\newtheorem{theorem}{Theorem}[section]

\newtheorem{lemma}[theorem]{Lemma}

\theoremstyle{definition}
\newtheorem{definition}[theorem]{Definition}

\theoremstyle{remark}
\newtheorem{remark}[theorem]{Remark}

\numberwithin{equation}{section} %% Equation numbering control.
\numberwithin{figure}{section}   %% Figure numbering control.

\usepackage[square,comma,numbers,sort&compress]{natbib}
\usepackage{url}

\newcounter{my_counter}
\begin{document}

\title[Cont. Data Assimilation for the 3D Primitive Eq.]{Continuous data assimilation for the 3D primitive equations of the ocean}

\date{\today}

% Group authors per affiliation:
\author{Yuan Pei}
\address[Yuan Pei]{Department of Mathematics, 
                University of Nebraska--Lincoln,
        Lincoln, NE 68588-0130, USA}
\email[Yuan Pei]{Yuan.Pei@wwu.edu}

\begin{abstract}
In this article, we show that the continuous data assimilation algorithm is valid for the 3D primitive equations of the ocean. Namely, the $L^2$ norm of the assimilated solution converge to that of the reference solution at an exponential rate in time. We also prove the global existence of strong solution to the assimilated system. 
\end{abstract}

\maketitle
\thispagestyle{empty}
%============================================================
\noindent

%
%\begin{keyword}
%Feedback control \sep
%Data assimilation\sep 
%Fluid dynamics \sep
%% Numerical analysis \sep
%% Mathematical computation \sep
%Primitive equations \sep 
%Regularity \sep
%Meteorology
%% Simulations \sep
%\MSC[2010]  
%35B35\sep
%35B40\sep
%35B41\sep
%35B45\sep
%35B65\sep
%35K40\sep
%35K51\sep
%35K61\sep
%35M13\sep
%35Q93\sep
%76D03\sep
%76D05\sep
%76D55\sep
%76D99
%\end{keyword}

%\end{frontmatter}

%\linenumbers

\section{Introduction}
\label{secInt}
\noindent
In this paper, we address the continuous data assimilation algorithm applied to the 3D primitive equations of the ocean. We first show that the data assimilation primitive system possesses a global in time strong solution. Then, we prove that the strong solution of the assimilated system converges to that of the reference system exponentially fast in both $L^2$ norm. Namely, the continuous data assimilation algorithm, which we elaborate in details below, is valid for the 3D primitive equations of the ocean. This work sets up the analytical foundation and provides the guideline for the application of the new and promising continuous data assimilation algorithm for the primitive equations of the ocean, which are widely used in geophysics and meteorology. 

We start with an introduction to the idea of data assimilation following \cite{FLT}. Simulations for dynamical systems, such as the primitive equations, require accurate initialization process in order to make reliable predictions. The initialization procedure depends on how acquired observations, such as measurements of temperature and velocity, are properly interpolated in space and time, to complete the information throughout the entire space-time domain, while maintaining the dynamical balance between these fields. Even if direct observations of these fields are available, typically they are not uniformly distributed in space or are very sparse. Moreover, the errors contained in the measurements and parameters of the model, combined with the highly nonlinear dynamics of the PDEs, makes basic interpolation unfitting for initial condition construction or selection, even for short term predictions. Meteorologists have designed a series of diagnostic experiments when creating precise initializing processes that minimize the loss of information from the measurements. For instance, the interpolated function in space and time must satisfy the conservation laws in the continuous model equations. Meteorologists also use various combinations of collected information about the state of the system to see which combination(s) will drive the system closer to the collected data. In this process, they may also use other forms of collected measurements such as, the average temperature of a small region obtained through image processing of satellite observation. One must know how to properly make use of this information, in particular when the collected data is not one of the evolving state variable in the dynamical system. In the context of meteorology and atmospheric physics, data assimilation algorithm when some state variable observations are not available as an input, has been studied in \cite{FLT,EB,HA,LAE} and the references therein, for some simplified numerical forecast models. Although nonlinear interactions between the scales of motion and the parameters in the system exist, it has been shown in several settings that if the dynamical model used as a source of {\it{a priori}} information captures the important properties of the system being modeled, then, one can identify the full state of the system knowing only coarse observation from partial data that is selected properly. 

In this work, we consider the continuous data assimilation introduced by Azouani, Olson, and Titi in \cite{AOT,AT} (see also \cite{CKT,HOT,OT} for early ideas in this direction). This approach, which we call AOT data assimilation or the linear AOT algorithm, is based on feedback control at the partial differential equation level. We now describe the idea of the AOT algorithm. Consider a dynamical system in a general form, 
\begin{equation}\label{ODE}
 \left\{
  \begin{aligned}
     \frac{d}{dt} u &= F(u),\\
    u(0)&=u_0.\\
  \end{aligned}
 \right.
\end{equation}
For example, this could represent a system of partial differential equations modeling fluid flow in the atmosphere or the ocean.  A central difficulty is that, even if one were able to solve the system exactly, the initial data $u_0$ is largely unknown.  For example, in a weather or climate simulation, the initial data may be measured at certain locations by weather stations, but the data at locations in between these stations remain unknown.  Therefore, one might not have access to the complete initial data $u_0$, but only to the observational measurements, which we denote by $I_{\delta}(u_0)$. (Here, $I_{\delta}$ is assumed to be a linear operator that can be taken, for example, to be an interpolation operator between grid points of maximal spacing $h$, or as an orthogonal projection onto Fourier modes no larger than $k\sim 1/h$.)
Moreover, the data from measurements may be streaming in moment by moment, so in fact, one often has the information $I_{\delta}(u)=I_{\delta}(u(t))$, for a range of times $t$.  Data assimilation is an approach that eliminates the need for complete initial data and also incorporates incoming data into simulations.  Classical approaches to data assimilation are typically based on the Kalman filter.  See, e.g., \cite{rD,eK,LSZ} and the references therein for more information about the Kalman filter.  In 2014, the AOT algorithm was introduced in \cite{AOT,AT}.  This new approach overcomes some of the drawbacks of the Kalman filter approach (see, e.g., \cite{BHLP} for further discussion). The approach has been the subject of much recent study in various contexts, see, e.g., \cite{ANLT,ATG,BOT,BM,FJT,FLT1,FLT2,FLT3,FMT,GOT,JMT,LP,MTT,MT}.

Specifically, the following system was proposed and studied in \cite{AOT,AT}: 
\begin{equation}\label{AOT}
 \left\{
  \begin{aligned} 
    \frac{d}{dt} \tilde{u} &= F(\tilde{u}) + \mu(I_{\delta}(u) - I_{\delta}(\tilde{u})),\\
    \tilde{u}(0) &= \tilde{u}_0.\\
  \end{aligned}
 \right.
\end{equation}
This system, used in conjunction with \eqref{ODE}, forms the AOT algorithm for data assimilation of system \eqref{ODE}.  In the case where the dynamical system \eqref{ODE} is the 2D Navier-Stokes equations, it was proven in \cite{AOT,AT} that, for any divergence-free initial data $\tilde{u}_0\in L^2$, 
$
 \|u(t) - \tilde{u}(t)\|_{L^2} \rightarrow 0,
$
exponentially in time. 
In particular, even without knowing the initial data $u_0$, the solution $u$ can be approximately reconstructed for large times.  We emphasize that, as noted in \cite{AOT}, the initial data for \eqref{AOT} can be any $L^2$ function, even $v_0=0$. Thus, no information about the initial data is required to reconstruct the solution asymptotically in time.

In the description below we denote the two-dimensional horizontal gradient, Laplacian, and divergence operators by $\nabla_{2}$, $\Delta_{2}$, and $div_{2}$, respectively. 
The 3D primitive equations of the ocean and atmosphere (PEs) on the domain $\Omega = \mathbb{R}^2/\mathbb{Z}^2 \times (-h, 0)$ with $h>0$ read:
\begin{equation}\label{PE}
 \left\{
  \begin{aligned}
   &\frac{\partial v}{\partial t} + (v\cdot\nabla_{2})v + w\frac{\partial v}{\partial z} + \nabla_{2} p + f\vec{k}\times v + L_{1}v = 0, 
      \\
   &\partial_{z}p + \theta = 0,
      \\
   &\nabla_{2}\cdot v + \partial_{z}w = 0,
      \\
   &\frac{\partial \theta}{\partial t} + v\cdot\nabla_{2} \theta + w\frac{\partial \theta}{\partial z} + L_{2}\theta = q,
      \\
   &\text{with } \quad v(x,0) =v_0(x), \quad \text{ and } \quad \theta(x, 0) = \theta_{0}(x), 
  \end{aligned}
 \right.
\end{equation}
where  $v = (v_1, v_2)$ and $w$ stand for the horizontal and vertical components of the three-dimensional velocity field $u = (v, w)$ and $\theta$ for the temperature and $p$ for the scalar pressure, all of which are unknowns of the equations. Also, $f$ is the Coriolos parameter and $q$ is the given heating source. The viscosity and the heat diffusion operators $L_{1}$ and $L_{2}$ are given by 
$$L_{1} = -\nu_{v}\Delta_{2} - \nu_{w}\frac{\partial^2}{\partial z^2} \quad \text{ and } \quad L_{2} = -\eta_{v}\Delta_{2} - \eta_{w}\frac{\partial^2}{\partial z^2},$$ 
where the positive constants $\nu_{v}$, $\nu_{w}$, and $\eta_{v}$, $\eta_{w}$ are horizontal and vertical viscosities and heat diffusivities, respectively. 
For simplicity, we assume that the initial data is sufficiently smooth (made more precise below) and satisfy $$div_{2}\int_{-h}^{0}v_0(x)\,dx_3 = 0$$ in the sense of distributions. Also, on top $z=0$ and $z=-h$ of the domain, we equip \eqref{PE} with Neumann boundary conditions 
$$\frac{\partial v}{\partial z} = 0, w = 0 \quad \text{ and } \quad \frac{\partial \theta}{\partial z} = 0.$$

The primitive equations of the ocean are considered to be the fundamental model for meteorology and climate prediction (c.f. \cite{jP}). The primitive equations are derived from the full compressible Navier-Stokes equations, which are too complicated in terms of application and contain phenomena such as shocks and sound waves that are not interesting from the perspective of geophysics. Thus, by applying the Boussinesq and the hydrostatic approximations to the full Navier-Stokes system, one obtains the primitive system consists of the horizontal momentum equations, the conservation of mass, and the thermodynamic equation--diffusion of temperature. 

The mathematical analysis of PE was initiated in \cite{LTW1, LTW2, LTW3} which set the analytical foundation for the system and established the global existence of weak solutions for $L^2$ initial data in the spirit of Leray. The $H^2$ regularity of the associated stationary linear problem was obtained in \cite{Z1, Z2}. This result has led to the local existence of strong solutions with initial data in $H^1$, which was established in \cite{BGMR1, BGMR2, BGMR3, GMR, HTZ}. The global existence of strong solutions with initial data in $H^1$ was proven by Cao and Titi in the classic paper \cite{CT1}. For other works on the primitive equations, we refer the readers to \cite{CT2,CLT1,CLT2,CINT,CLT3,CLT4,GH,GHZ,H,K,LT,KPRZ,KZ1,KZ2,KZ3,P,STT,SV} as well as \cite{KTVZ1,KTVZ2,MW,R,RTT1,RTT2} for the inviscid case. 

From the point of view of data assimilation, the primitive equations are ideal in applying the algorithm since the solution to the system has chaotic large-time behavior; and it also has a finite-dimensional global attractor (see \cite{C,J1,JT} and the reference therein) making PE an excellent candidate for studying large-time behavior. Therefore, in view of the AOT data assimilation, we have:  
\begin{equation}\label{PE_DA}
 \left\{
  \begin{aligned}
   &\frac{\partial \tilde{v}}{\partial t} 
   + 
   (\tilde{v}\cdot\nabla_{2})\tilde{v} 
   + 
   \tilde{w}\frac{\partial \tilde{v}}{\partial z} 
   + 
   \nabla_{2} \tilde{p} 
   + 
   f\vec{k}\times \tilde{v} + L_{1}\tilde{v} 
   =
   \mu_{u}(I_{\delta}(v) - I_{\delta}(\tilde{v})), 
      \\
   &\partial_{z}\tilde{p} + \tilde{\theta} = 0,
      \\
   &\nabla_{2}\cdot \tilde{v} + \partial_{z}\tilde{w} = 0,
      \\
   &\frac{\partial \tilde{\theta}}{\partial t} + \tilde{v}\cdot\nabla_{2} \tilde{\theta} + \tilde{w}\frac{\partial \tilde{\theta}}{\partial z} + L_{2}\tilde{\theta} 
   = 
   q + \mu_{\theta}(I_{\delta}(\theta) - I_{\delta}(\tilde{\theta})),
      \\
   &\text{with } \quad \tilde{v}(x,0) = 0, \quad \text{ and } \quad \tilde{\theta}(x, 0) = 0, 
  \end{aligned}
 \right.
\end{equation}
where $1/\delta$ is the spatial resolution and the operator $I_{\delta}$ satisfies the following inequality:
\begin{align}
   \label{C1}
   \Vert w-I_\delta(w)\Vert_{L^2} &\leq c_1\delta\Vert\nabla w\Vert_{L^2} \quad\text{for all}\quad w\in H^1,
   \\
   \label{C2}
   \Vert I_{\delta}(w)\Vert_{L^2}& \leq C\Vert w\Vert_{L^2} \quad\text{for all}\quad w\in L^2.
\end{align}    
In particular, one may take $I_\delta$ to be orthogonal projection onto the first $c/\delta$ Fourier modes, for some constant $c$.  
\begin{remark}
Other physically relevant choices of $I_\delta$, such as a nodal interpolation operator, have been considered in the case of the linear AOT algorithm (see, e.g., \cite{GOT}). 
\end{remark}

The principal aim of this article is double-fold. First, we show the existence of global strong solution to the assimilation system, which is necessary for the convergence results. For the proof, we follow the idea in the classic paper \cite{CT1} and provide the {\it{a priori}} estimates for the $L^2$, $L^6$, and $H^1$ norms of the solutions. The main difficult comes from the lack of the momentum equation for the third component of the velocity. We use the anisotropic estimates as well as the divergence-free condition to bound the nonlinear terms. Secondly, to the best knowledge of the author, we present the first continuous data assimilation result for the full primitive equations of the ocean in 3D. We prove that with a sufficiently refined resolution ($\delta$ small) and large enough nudging coefficients ($\mu_{u}$ and $\mu_{\theta}$), the assimilated solution approaches the reference solution at an exponential rate in the sense of $L^2$ norm. 

\section{Preliminaries and Summary of Main Results}
\label{Pre}
\noindent
In this paper, we frequently use the $L^2$ norm and $H^1$ semi-norms, defined by
\begin{align*}
 \|u\|_{L^2}^2 = \int_{\Omega}|u(x)|^2\,dx,
 \qquad
 \|u\|_{H^1} = \|\nabla u\|_{L^2}.
\end{align*}

All through this paper we denote 
$\partial_{j} = \partial / \partial x_{j}$, 
$\partial_{jj} = \partial^{2}/\partial x_{j}^2$, 
$\partial_{t} = \partial /\partial t$, 
and the horizontal gradient $\nabla_{2} = (\partial_1, \partial_2)$ 
and horizontal Laplacian $\Delta_{2} = \partial_{11} + \partial_{22}$. 
Also, we denote the usual Lebesgue and Sobolev spaces by $L_{x}^{p}$ and $H_{x}^{s}\equiv W_{x}^{s, 2}$, respectively, 
and let $\mathcal{F}$ be the set of all trigonometric polynomial over $\mathbb{T}^3$ 
and define the space of zero-average smooth functions 
$$\mathcal{V} := \left\{ \phi\in\mathcal{F}\times\mathcal{F}: \nabla_2\cdot\int_{-h}^{0}\phi\,dz = 0, \left.{\frac{\partial \phi}{\partial z}}\right|_{z=-h} = \left.{\frac{\partial \phi}{\partial z}}\right|_{z=-0} = 0, \right\}.$$
Also, we define $H$ and $V$ to be the closures of $\mathcal{V}$ in $L_{x}^2$ and $H_{x}^1$, 
with inner products 
$$(u, v) = \sum_{i=1}^3\int_{\mathbb{T}^3}u_{i}v_{i}\,dx \text{ \,\,and\,\,  } ((u, v)) = (\nabla u, \nabla v) = \sum_{i, j=1}^3\int_{\Omega}\partial_{j}u_{i}\partial_{j}v_{i}\,dx,$$
respectively, associated with the norms $|u|=(u, u)^{1/2}$ and $\Vert u \Vert=(( u, u))^{1/2}$.
The latter is indeed a norm due to Poincar\'e inequality $$\lambda_1\Vert\phi\Vert_{L_{x}^2}^2 \leq \Vert\nabla\phi\Vert_{L_{x}^2}^2$$ 
for $\phi\in V$ where $\lambda_1$ is the first eigenvalue of the Stokes operator $A := -P_{\sigma}\Delta$. 
Similarly we define $$\mathcal{V}_{\theta} := \left\{ \phi\in\mathcal{F}: \left.{\frac{\partial \phi}{\partial z}}\right|_{z=-h} = \left.{\frac{\partial \phi}{\partial z}}\right|_{z=-0} = 0, \int_{\Omega} \phi = 0\right\},$$ 
and denote the closure of $\mathcal{V}_{\theta}$ in $L^2$ and $H^1$ by $H_{\theta}$ and $V_{\theta}$, for $\theta$. 
We refer to readers to \cite{CF, T} for more standard details on Navier-Stokes, which will be used in this paper. 

The following Gagliardo-Nirenberg-Sobolev type inequality is frequently used in our estimates. 
\begin{lemma}
\label{L1}
Assume $1 \leq q, r \leq \infty$, and $0<\gamma<1$.  
For $v\in L_{x}^q(\mathbb{T}^{n})$, such that  $\partial^\alpha v\in L_{x}^{r} (\mathbb{T}^{n})$, for $|\alpha|=m$, then 
\begin{align}\label{PT1}
\Vert\partial_{s}v\Vert_{L^{p}} \leq C\Vert\partial^{\alpha}v\Vert_{L^{r}}^{\gamma}\Vert v\Vert_{L^{q}}^{1-\gamma},
\quad\text{where}\quad
\frac{1}{p} - \frac{s}{n} = \left(\frac{1}{r} - \frac{m}{n}\right) \gamma+ \frac{1}{q}(1-\gamma).
\end{align}
\end{lemma}

We list some special cases below, 
\begin{lemma}
\label{L2}
For $\phi \in H^1(\mathbb{T}^2=\mathbb{R}^2/\mathbb{Z}^2)$, 
we have 
$$\Vert\phi\Vert_{L^4} \leq C\Vert\phi\Vert_{L^2}^{1/2}\Vert\phi\Vert_{H^1}^{1/2}, $$
$$\Vert\phi\Vert_{L^8} \leq C\Vert\phi\Vert_{L^6}^{3/4}\Vert\phi\Vert_{H^1}^{1/4}; $$
for $u \in H^1(\Omega)$, we have 
$$\Vert u\Vert_{L^3} \leq C\Vert u\Vert_{L^2}^{1/2}\Vert u\Vert_{H^1}^{1/2}, $$
$$\Vert u\Vert_{L^6} \leq C\Vert u\Vert_{H^1}, $$
and the constant $C$ may vary but is independent of $\phi$ and $u$, respectively. 
\end{lemma}

The next lemma about anisotropic estimates is critical to out treatment of the nonlinear terms.  
\begin{lemma}
\label{L3}
For $v \in V\times V$, $\phi \in H^1$, and $\psi \in L^2$, 
there exists constant $C$ independent of $v, \phi, \psi$ such that  
\begin{equation*}
     \int\left(\int_{-h}^{z} \nabla_{2}\cdot v(\cdot, \zeta)\,d\zeta\right)\,\phi\,\psi 
     \leq
     C\Vert\nabla v\Vert_{L^2}^{1/2} \Vert\nabla v\Vert_{H^1}^{1/2} \Vert \phi\Vert_{L^2}^{1/2} \Vert\phi\Vert_{H^1}^{1/2} \Vert\psi\Vert_{L^2}^{1/2}. 
\end{equation*}
\end{lemma}

Now we state our main results of this paper. 
The first theorem concerns the global existence of solution to system \eqref{PE_DA}. 
\begin{theorem}
\label{T0} 
For $\tilde{v}_0 \in V$, $\tilde{\theta}_0 \in V_{\theta}$, and $q \in H^1$, with $\mu_{u}, \mu_{\theta} > 0$, system \eqref{PE_DA} possesses a unique global-in-time strong solution. 
\end{theorem}

The next theorem provides the exponential-in-time convergence in $L^2$-norm of the continuous data assimilation for the 3D primitive equations of the ocean. 
\begin{theorem}
\label{T4}
Suppose $v_0 \in V, \theta_0 \in V_{\theta}, q \in H^1$ and $\tilde{v}_0 = 0$, $\tilde{\theta}_0 = 0$. Assume condition \eqref{C1} holds, with 
$$\mu_{u}, \mu_{\theta} > 2(C + C(\Vert v\Vert_{H^2}^4 + \Vert \theta\Vert_{H^2}^4))$$ 
where $C$ depends on $\lambda_1, \nu_{v}, \nu_{w}, \eta_{v}, \eta_{w}$, the domain $\Omega$, and the reference solution, and with 
$$0 < \delta < \frac{1}{\sqrt{2}c_1}\frac{\sqrt{\mu_{min}}}{\mu_{max}\sqrt{\nu_{v}^{-1}+\nu_{w}^{-1}+\eta_{v}^{-1}+\eta_{w}^{-1}}},$$ 
we have exponential-in-time convergence in $L^2$ norm of the solution of \eqref{PE_DA} to that of \eqref{PE}, i.e., $\Vert v(t) - \tilde{v}(t)\Vert_{L^2} \to 0$ and $\Vert \theta(t) - \tilde{\theta}(t)\Vert_{L^2} \to 0$ at an exponential rate as $t \to \infty$. 
\end{theorem}

Recall the definition of the strong solution to \eqref{PE} as follows, 
\begin{definition}
\label{SS}
For $v_0 \in V$, $\theta_0 \in V_{\theta}$, and $T>0$, we say that $(v, \theta)$ is the strong solution of \eqref{PE} on the time interval $[0, T]$ if it satisfies \eqref{PE} in the weak sense and also we have 
\begin{align*}
     v & \in C([0, T]; V) \cap L^{2}([0, T], H^2(\Omega)\cap V),
     \\
     \theta & \in C([0, T]; V_{\theta}) \cap L^{2}([0, T], H^2(\Omega)\cap V_{\theta}),
     \\
     \frac{d v}{d t} & \in L^{1}([0, T]; H),
     \\
     \frac{d \theta}{d t} & \in L^{1}([0, T]; H_{\theta}). 
\end{align*}
\end{definition}
In order to prove our main theorems, we need the following theorem summarized from the results in \cite{CT1} on the global well-posedness of system \eqref{PE}.
\begin{theorem}
\label{T1}
Let $Q\in H^1$, $v_0 \in V$, $\theta_0\in V_{\theta}$ and $T>0$. 
Then, there exists a unique strong solution of \eqref{PE} on $[0, T]$ which continuous depends on the initial data. 
\end{theorem}

Further, we summarize the following theorems from \cite{C,J1,JT} regarding the uniform bounds of the $L^2$, $H^1$, and $H^2$ norms of the solution and the existence of finite dimensional global attractors of system \eqref{PE}. See also \cite{J2} for a new proof. 
\begin{theorem}
\label{T2}
For $v_0 \in V$, $\theta_0 \in V_{\theta}$, and $Q \in H^1$, 
the strong solution $(v, \theta)$ satisfies 
\begin{align*}
     v & \in L^{\infty}([0, \infty); V) \cap L^{\infty}([0, \infty); L^6\cap H)
     \\
     \theta & \in L^{\infty}([0, \infty); V_{\theta}) \cap L^{\infty}([0, \infty); L^6\cap H_{\theta}).
\end{align*}
Moreover, there exists a bounded absorbing ball for all solutions $(v, \theta)$ of \eqref{PE} in $V\times V_{\theta}$. 
\end{theorem}

\begin{theorem}
\label{T3}
With the same assumption as in Theorem~\ref{T2}, 
for the strong solution $(v, \theta)$ of \eqref{PE}, 
we have 
\begin{align*}
     v & \in L^{\infty}([0, \infty); H^2\cap V),
     \\
     \theta & \in L^{\infty}([0, \infty); H^2\cap V_{\theta}),
     \\
     \frac{\partial v}{\partial z} & \in L^{2}([0, \infty); V),
     \\
     \frac{\partial \theta}{\partial z} & \in L^{2}([0, \infty); V_{\theta}).
\end{align*}
Furthermore, there exists a bounded absorbing ball for all solutions of \eqref{PE} in $(H^2\cap V)\times (H^2\cap V_{\theta})$ and system \eqref{PE} possesses a compact and connected global attractor with finite fractal and Hausdorff dimensions. 
\end{theorem}

For the sake of completeness, we state a special case of the generalized Gr\"onwall's inequality used in \cite{FLT3} and  \cite{daJT}. 
\begin{lemma}
\label{L4}
Suppose that $Y(t)$ is a locally integrable and absolutely continuous function that satisfies the following: 
$$\frac{d Y}{d t} + \alpha(t) Y \leq \beta(t), \quad\text{ a.e. on } (0, \infty), $$
such that 
$$\liminf_{t \to \infty} \int_{t}^{t+\tau} \alpha(s)\,ds \geq \gamma, \quad\quad\quad \limsup_{t \to \infty} \int_{t}^{t+\tau} \alpha^{-}(s)\,ds < \infty, $$
and 
$$\lim_{t \to \infty} \int_{t}^{t+\tau} \beta^{+}(s)\,ds = 0, $$
for fixed $\tau > 0$, and $\gamma > 0$, 
where 
$\alpha^{-} = \max\{-\alpha, 0\}$ 
and 
$\beta^{+} = \max\{\beta, 0\}$. 
Then, $Y(0) \to 0$ at an exponential rate as $t \to \infty$.
\end{lemma}

\section{Existence of solution to the assimilated system}
\label{Exist_Uniq}
In this section, we provide the proof for the existence of solutions to system \eqref{PE_DA}. 
We provide the {\it{a priori}} estimates only and one can obtain the solution via standard Galerkin approximation henceforth, which we omit for simplicity. 
\begin{proof}
{\smallskip\noindent {\em Proof of Theorem~\ref{T0}.}} 
\subsection{$L^2$-estimates}
First we multiply $\tilde{v}$ and $\tilde{\theta}$ to the corresponding equations in system \eqref{PE_DA}, respectively, integrate over $\Omega$, and get 
\begin{equation}\label{L2est}
 \left\{
   \begin{aligned}
   &\frac{1}{2}\frac{d}{d t}\Vert \tilde{v}\Vert_{L^2}^2 
   +
   \nu_{v}\Vert\nabla_{2}\tilde{v}\Vert_{L^2}^2 
   + 
   \nu_{w}\Vert\partial_{z}\tilde{v}\Vert_{L^2}^2
   = 
   -\int_{\Omega}(\tilde{v}\cdot\nabla_{2})\tilde{v}\cdot \tilde{v} 
   +
   \int_{\Omega}\tilde{w}\partial_{z}\tilde{v}\cdot \tilde{v} 
   \\&\qquad
   + 
   \int_{\Omega}(\nabla_{2}\int_{-h}^{z}\tilde{\theta}\,d\zeta)\cdot \tilde{v} 
   +
   \mu_{u}\int_{\Omega}(I_{\delta}(v) - I_{\delta}(\tilde{v}))\cdot \tilde{v}, 
   \\
   &\frac{1}{2}\frac{d}{d t}\Vert \tilde{\theta}\Vert_{L^2}^2 
   +
   \eta_{v}\Vert\nabla_{2}\tilde{\theta}\Vert_{L^2}^2 + \eta_{w}\Vert\partial_{z}\tilde{\theta}\Vert_{L^2}^2 
   +
   \Vert\tilde{\theta}(z=0)\Vert_{L^2}^2
   = 
   -\int_{\Omega}\tilde{\theta}\left(\tilde{v}\cdot\nabla_{2} \tilde{\theta}\right)
    \\&\qquad
   + 
   \int_{\Omega}\left(\int_{-h}^{z}\nabla_{2}\cdot \tilde{v}\,d\zeta\right)\tilde{\theta}\partial_{z}\tilde{\theta} 
   +
   \int_{\Omega}q\tilde{\theta}
   + 
   \mu_{\theta}\int_{\Omega}(I_{\delta}(\theta) - I_{\delta}(\tilde{\theta}))\tilde{\theta},
  \end{aligned}
 \right.
\end{equation}
where we used 
$$\tilde{p} = -\int_{-h}^{z} \tilde\Theta\,d\zeta + \tilde{p}_0, $$
and without loss of generality, we assume $\tilde{p}_0 = 0$. 
Due to the divergence-free condition, the first two terms on the right side of the first equation above add up to zero. 
The third term is bounded by 
\begin{align*}
     &\int_{\Omega}|\tilde{\theta}| |\nabla\tilde{v}| 
     \leq
     C\Vert\tilde{\theta}\Vert_{L^2} \Vert\nabla\tilde{v}\Vert_{L^2}, 
\end{align*}
while for the assimilation term, we observe that 
\begin{align*}
     &\Vert I_{\delta}(v)\Vert_{L^2}
     \leq
     \Vert I_{\delta}(v) - v\Vert_{L^2}
     +
     \Vert v\Vert_{L^2}
     \leq
     (c_1\delta + \lambda_1^{-1/2})\Vert\nabla v\Vert_{L^2}, 
\end{align*}
where we used condition \eqref{C1} and Poincar\'e's inequality. 
Thus, the last term is bounded by 
\begin{align*}
     &\mu_{u}\int_{\Omega}I_{\delta}(v)\tilde{v} 
     -
     \mu_{u}\int_{\Omega}I_{\delta}(\tilde{v})\tilde{v}
     \\&
     \leq
     \mu_{u}\Vert I_{\delta}(v)\Vert_{L^2} \Vert\tilde{v}\Vert_{L^2}
     -
     \mu_{u}\int_{\Omega}(I_{\delta}(\tilde{v}) - \tilde{v})\tilde{v}
     -
     \mu_{u}\Vert\tilde{v}\Vert_{L^2}^2
     \\&
     \leq
     C\mu_{u}\Vert\nabla v\Vert_{L^2} \Vert\tilde{v}\Vert_{L^2}
     +
     c_1\delta\mu_{u}\Vert\nabla\tilde{v}\Vert_{L^2} \Vert\tilde{v}\Vert_{L^2} 
     -
     \mu_{u}\Vert\tilde{v}\Vert_{L^2}^2
     \\&
     \leq
     C\Vert\nabla v\Vert_{L^2}^2
     +
     \frac{\mu_{u}}{2}\Vert\tilde{v}\Vert_{L^2}^2
     +
     \frac{\mu_{u}}{2}c_1^2\delta^2\Vert\nabla\tilde{v}\Vert_{L^2}^2
     +
     \frac{\mu_{u}}{2}\Vert\tilde{v}\Vert_{L^2}^2
     -
     \mu_{u}\Vert\tilde{v}\Vert_{L^2}^2.
\end{align*} 
Now choose $\delta$ small enough so that 
$$\frac{\mu_{u}}{2}c_1^2\delta^2 < \frac{\min{\{\nu_{v}, \nu_{w}\}}}{4}, $$
we obtain 
\begin{align*}
     &\frac{1}{2}\frac{d}{d t}\Vert \tilde{v}\Vert_{L^2}^2 
     +
     \frac{\nu_{v}}{2}\Vert\nabla_{2}\tilde{v}\Vert_{L^2}^2 + \frac{\nu_{w}}{2}\Vert\partial_{z}\tilde{v}\Vert_{L^2}^2
     \leq
     C\Vert\tilde{\theta}\Vert_{L^2}^2
     +
     C\Vert\nabla v\Vert_{L^2}^2, 
\end{align*} 
where we used Cauchy-Schwarz inequality. 
Regarding the estimates of the right side of $\tilde{\theta}$-equation, 
we also observe that divergence-free condition implies that the first two terms sum up to zero, 
while the third term is bounded by 
\begin{align*}
     &\int_{\Omega}|q| |\tilde{\theta}|
     \leq
     \Vert q\Vert_{L^2} \Vert\tilde{\theta}\Vert_{L^2}. 
\end{align*}
As for the last term, similar to the case in $\tilde{v}$, we bound it by 
\begin{align*}
     &\mu_{\theta}\int_{\Omega}I_{\delta}(\theta)\tilde{\theta} 
     -
     \mu_{\theta}\int_{\Omega}I_{\delta}(\tilde{\theta})\tilde{\theta}
     \\&
     \leq
     \mu_{\theta}\Vert I_{\delta}(\theta)\Vert_{L^2} \Vert\tilde{\theta}\Vert_{L^2}
     -
     \mu_{\theta}\int_{\Omega}(I_{\delta}(\tilde{\theta}) - \tilde{\theta})\tilde{\theta}
     -
     \mu_{\theta}\Vert\tilde{\theta}\Vert_{L^2}^2
     \\&
     \leq
     C\mu_{\theta}\Vert\nabla\theta\Vert_{L^2} \Vert\tilde{\theta}\Vert_{L^2}
     +
     c_1\delta\mu_{\theta}\Vert\nabla\tilde{\theta}\Vert_{L^2} \Vert\tilde{\theta}\Vert_{L^2} 
     -
     \mu_{\theta}\Vert\tilde{\theta}\Vert_{L^2}^2
     \\&
     \leq
     C\Vert\nabla\theta\Vert_{L^2}^2
     +
     \frac{\mu_{\theta}}{2}\Vert\tilde{\theta}\Vert_{L^2}^2
     +
     \frac{\mu_{\theta}}{2}c_1^2\delta^2\Vert\nabla\tilde{\theta}\Vert_{L^2}^2
     +
     \frac{\mu_{\theta}}{2}\Vert\tilde{\theta}\Vert_{L^2}^2
     -
     \mu_{\theta}\Vert\tilde{\theta}\Vert_{L^2}^2.
\end{align*}
By choosing small enough $\delta$ such that 
$$\frac{\mu_{\theta}}{2}c_1^2\delta^2 < \frac{\min{\{\eta_{v}, \eta_{w}\}}}{4}, $$
we get 
\begin{align*}
     &\frac{1}{2}\frac{d}{d t}\Vert \tilde{\theta}\Vert_{L^2}^2 
     +
     \frac{\eta_{v}}{2}\Vert\nabla_{2}\tilde{\theta}\Vert_{L^2}^2 + \frac{\eta_{w}}{2}\Vert\partial_{z}\tilde{\theta}\Vert_{L^2}^2 
     +
     \Vert\tilde{\theta}(z=0)\Vert_{L^2}^2
     \\&
     \leq
     C\Vert\nabla\theta\Vert_{L^2}^2 
     +
     \Vert q\Vert_{L^2} \Vert\tilde{\theta}\Vert_{L^2}
     \\&
     \leq
     C\Vert\nabla\theta\Vert_{L^2}^2
     +
     \Vert q\Vert_{L^2}\sqrt{2}h\Vert\partial_{z}\tilde{\theta}\Vert_{L^2}
     +
     \Vert q\Vert_{L^2}\sqrt{2h}\Vert\tilde{\theta}(z=0)\Vert_{L^2}
     \\&
     \leq
     C\Vert\nabla\theta\Vert_{L^2}^2
     +
     \frac{C_{h}}{\eta_{w}}\Vert q\Vert_{L^2}^2
     +
     \frac{\eta_{w}}{4}\Vert\partial_{z}\tilde{\theta}\Vert_{L^2}^2
     +
     \frac{1}{2}\Vert\tilde{\theta}(z=0)\Vert_{L^2}^2,
\end{align*}
where we used the fact that 
$$\Vert\tilde{\theta}\Vert_{L^2}^2 \leq 2h^2\Vert\partial_{z}\tilde{\theta}\Vert_{L^2}^2 + 2h\Vert\tilde{\theta}(z=0)\Vert_{L^2}^2. $$
Therefore, by Theorem~\ref{T2} and Gr\"onwall's inequality, we conclude that 
$$\Vert\tilde{\theta}\Vert_{L^2} \leq K_1,$$ 
where $K_1$ depends only on $\Vert\theta\Vert_{H^1}$ and $\Vert q\Vert_{L^2}$. 
Whence, we obtain 
$$\Vert\tilde{v}\Vert_{L^2} \leq K_2,$$ 
and here $K_2$ depends on $K_1$ and $\Vert v\Vert_{H^1}$. 
Furthermore, we have the bounds 
\begin{align*}
     &\int_{t=0}^{t=T}\eta_{v}\Vert\nabla_{2}\tilde{\theta}\Vert_{L^2}^2 + \eta_{w}\Vert\partial_{z}\tilde{\theta}\Vert_{L^2}^2 
     +
     \Vert\tilde{\theta}(z=0)\Vert_{L^2}^2 \,dt
     \leq 
     2K_{1}T
\end{align*}
and
\begin{align*}
     &\int_{t=0}^{t=T}\nu_{v}\Vert\nabla_{2}\tilde{v}\Vert_{L^2}^2 
     + 
     \nu_{w}\Vert\partial_{z}\tilde{v}\Vert_{L^2}^2 \,dt
     \leq 
     2K_{2}T
\end{align*}
for any $T>0$.

\subsection{$L^6$-estimates}
We begin by multiplying $|\tilde{v}|^{4}\tilde{v}$ and $|\tilde{\theta}|^{4}\tilde{\theta}$ to the equations of $\tilde{v}$ and $\tilde{\theta}$, respectively, integrate over $\Omega$, and get 
\begin{equation}\label{L6est}
 \left\{
   \begin{aligned}
     &\frac{1}{6}\frac{d}{dt}\Vert \tilde{v}\Vert_{L^6}^6
     +
     \frac{9\nu_{v}}{5}\Vert\nabla_2(|\tilde{v}|^3)\Vert_{L^2}^2
     +
     \frac{9\nu_{w}}{5}\Vert\partial_{z}(|\tilde{v}|^3)\Vert_{L^2}^2
     \\&
     =
     -
     \int_{\Omega}(\tilde{v}\cdot\nabla_2)\tilde{v}\cdot \tilde{v}|\tilde{v}|^4 
     -
     \int_{\Omega}\tilde{w}\partial_{z}\tilde{v}\cdot \tilde{v}|\tilde{v}|^4 
     + 
     \int_{\Omega}(\nabla_{2}\int_{-h}^{z}\tilde{\theta}\,d\zeta)\cdot \tilde{v}|\tilde{v}|^4 
     \\&\quad
     +
     \mu_{u}\int_{\Omega}(I_{\delta}(v) - I_{\delta}(\tilde{v}))\tilde{v}|\tilde{v}|^4, 
     \\&
     \frac{1}{6}\frac{d}{d t}\Vert \tilde{\theta}\Vert_{L^6}^6 
     +
     \frac{9\eta_{v}}{5}\Vert\nabla_{2}(|\tilde{\theta}|^3)\Vert_{L^2}^2
     +
     \frac{9\eta_{w}}{5}\Vert\partial_{z}(|\tilde{\theta}|^3)\Vert_{L^2}^2 
     +
     \Vert\tilde{\theta}(z=0)\Vert_{L^6}^6
     \\&
     = 
     - 
     \int_{\Omega}(\tilde{v}\cdot\nabla_{2}\tilde{\theta})|\tilde{\theta}|^{4}\tilde{\theta}
     - 
     \int_{\Omega}\tilde{w}\partial_{z}\tilde{\theta}|\tilde{\theta}|^{4}\tilde{\theta}
     +
     \int_{\Omega}q|\tilde{\theta}|^4\tilde{\theta}
      \\&\quad
     + 
     \mu_{\theta}\int_{\Omega}(I_{\delta}(\theta) - I_{\delta}(\tilde{\theta}))|\tilde{\theta}|^{4}\tilde{\theta}, 
  \end{aligned}
 \right.
\end{equation}
where we used the fact that $f\vec{k}\times\tilde{v}|\tilde{v}|^4 = 0$. 
Thanks to the divergence-free condition $\partial_{z}W = -\nabla_{2}\cdot \tilde{v}$, 
the first and second terms on the right side of both equations add up to zero, respectively.
Then, we estimate the remaining terms on the right side of the above two equations. 
For the third term on the right side of $\tilde{v}$-equation, we use H\"older's inequality and obtain    
\begin{align*}
     &\int_{\Omega}(\nabla_{2}\int_{-h}^{z}\tilde{\theta}\,d\zeta)\cdot \tilde{v}|\tilde{v}|^4
     =
     -\sum_{i=1}^2\int_{\Omega}(\int_{-h}^{z}\tilde{\theta}\,d\zeta)\cdot\partial_{i}\tilde{v}|\tilde{v}|^4
     \\&
     \leq
     C\Vert\tilde{\theta}\Vert_{L^6} \Vert\frac{5}{3}|\tilde{v}|^2\nabla(|\tilde{v}|^{2}\tilde{v})\Vert_{L^{6/5}}
     \leq
     C\Vert\tilde{\theta}\Vert_{L^6} \Vert|\tilde{v}|^2\Vert_{L^3} \Vert\nabla(|\tilde{v}|^3)\Vert_{L^2}
     \\&
     \leq
     C\left(\frac{1}{\nu_{v}}+\frac{1}{\nu_{w}}\right)\Vert\tilde{\theta}\Vert_{L^6}^6 
     +
     C\left(\frac{1}{\nu_{v}}+\frac{1}{\nu_{w}}\right)\Vert \tilde{v}\Vert_{L^6}^6
     \\&\quad
     +
     \frac{\nu_{v}}{8}\Vert\nabla_2(|\tilde{v}|^3)\Vert_{L^2}^2
     +
     \frac{\nu_{w}}{8}\Vert\partial_{z}(|\tilde{v}|^3)\Vert_{L^2}^2.
\end{align*}
The fourth term is bounded by  
\begin{align*}
     &C\mu_{u}\left( \Vert I_{\delta}(v) \Vert_{L^2} + \Vert I_{\delta}(\tilde{v}) \Vert_{L^2} \right) \Vert|\tilde{v}|^2\Vert_{L^3} \Vert(|\tilde{v}|^3)\Vert_{L^6}
     \\&     
     =
     C\mu_{u}\left( \Vert I_{\delta}(v) \Vert_{L^2} + \Vert I_{\delta}(\tilde{v}) \Vert_{L^2} \right) \Vert\tilde{v}\Vert_{L^6}^2 \Vert\nabla(|\tilde{v}|^3)\Vert_{L^2}
     \\&
     \leq
     C\left(\frac{1}{\nu_{v}}+\frac{1}{\nu_{w}}\right)\mu_{u}^6\left(\Vert v\Vert_{L^2}^6 + \Vert \tilde{v}\Vert_{L^2}^6\right)
     +
     C\Vert\tilde{v}\Vert_{L^6}^6
     \\&\quad
     +
     \frac{\nu_{v}}{8}\Vert\nabla_{2}(|\tilde{v}|^3)\Vert_{L^2}^2
     +
     \frac{\nu_{w}}{8}\Vert\partial_{z}(|\tilde{v}|^3)\Vert_{L^2}^2, 
\end{align*}
where we used condition \eqref{C2} for the operator $I_{\delta}$ and Cauchy-Schwarz inequality. 
Next, we estimate the third and fourth terms on the right side of the $\tilde{\theta}$-equation. 
By H\"older's inequality we first bound the third term by 
\begin{align*}
     &\Vert q\Vert_{L^6} \Vert\tilde{\theta}^5\Vert_{L^{6/5}}
     =
     \Vert q\Vert_{L^6} \Vert\tilde{\theta}\Vert_{L^6}^5
     \leq
     C\Vert q\Vert_{L^6}^6
     +
     C\Vert\tilde{\theta}\Vert_{L^6}^6
     \leq
     C\Vert q\Vert_{H^1}^6
     +
     C\Vert\tilde{\theta}\Vert_{L^6}^6
\end{align*}
while regarding the last term, we follow similarly the estimates for $\tilde{v}$ and bound it by  
\begin{align*}
     &C\mu_{\theta}\left( \Vert I_{\delta}(\theta) \Vert_{L^2} + \Vert I_{\delta}(\tilde{\theta}) \Vert_{L^2} \right) \Vert|\tilde{\theta}|^2\Vert_{L^3} \Vert(|\tilde{\theta}|^3)\Vert_{L^6}
     \\&     
     =
     C\mu_{\theta}\left( \Vert I_{\delta}(\theta) \Vert_{L^2} + \Vert I_{\delta}(\tilde{\theta}) \Vert_{L^2} \right) \Vert\tilde{\theta}\Vert_{L^6}^2 \Vert\nabla(|\tilde{\theta}|^3)\Vert_{L^2}
     \\&
     \leq
     C\left(\frac{1}{\eta_{v}}+\frac{1}{\eta_{w}}\right)\mu_{\theta}^6\left(\Vert \theta\Vert_{L^2}^6 + \Vert \tilde{\theta}\Vert_{L^2}^6\right)
     +
     C\Vert\tilde{\theta}\Vert_{L^6}^6
     \\&\quad
     +
     \frac{\eta_{v}}{8}\Vert\nabla_{2}(|\tilde{\theta}|^3)\Vert_{L^2}^2
     +
     \frac{\eta_{w}}{8}\Vert\partial_{z}(|\tilde{\theta}|^3)\Vert_{L^2}^2. 
\end{align*}
Thanks to Gr\"onwall's inequality, by summing up the above estimates and by Theorem~\ref{T1} and Theorem~\ref{T2}, we get 
$$\Vert\tilde{\theta}\Vert_{L^6} \leq K_3 \quad\text{ and }\quad \Vert\tilde{v}\Vert_{L^6} \leq K_4$$
where $K_3$ depends on $\Vert\theta\Vert_{H^1}$, $\Vert q\Vert_{H^1}$, and $c_1, \delta, \mu_{\theta}, \lambda_1$, while $K_4$ depends on $\Vert v\Vert_{H^1}$ and $K_3$. 

\subsection{$H^1$-estimates}
In order to obtain the $H^1$-bounds on $\tilde{v}$, 
we first need to estimate $\Vert\partial_{z}\tilde{v}\Vert_{L^2}$ and $\Vert\nabla\partial_{z}\tilde{v}\Vert_{L^2}$. 
\subsubsection{$\Vert\partial_{z}\tilde{v}\Vert_{L^2}$ estimates}
Multiply the first equation in \eqref{PE_DA} by $-\partial_{zz}^2\tilde{v}$, integrate over $\Omega$, and we get 
\begin{align*}
     &\frac{1}{2}\frac{d}{dt}\Vert\partial_{z}\tilde{v}\Vert_{L^2}^2
     +
     \nu_{v}\Vert\nabla_{2}\partial_{z}\tilde{v}\Vert_{L^2}^2
     +
     \nu_{w}\Vert\partial_{zz}^2\tilde{v}\Vert_{L^2}^2
     \\&
     =
     \int_{\Omega}(\tilde{v}\cdot\nabla_{2})\tilde{v}\partial_{zz}^2\tilde{v}
     +
     \int_{\Omega}\tilde{w}\partial_{z}\tilde{v}\partial_{zz}^2\tilde{v}
     \\&\quad
     -\int_{\Omega}(\nabla_{2}\int_{-h}^{z}\tilde{\theta}\,d\zeta)\cdot \partial_{zz}^2\tilde{v}
     -
     \mu_{u}\int_{\Omega}(I_{\delta}(v) - I_{\delta}(\tilde{v}))\partial_{zz}^2\tilde{v},
\end{align*}
where we used the equation $\partial_{z}\tilde{p} = -\tilde{\theta}$ and the fact that 
\begin{align*}
     &-\int_{\Omega} (f\vec{k}\times \tilde{v})\cdot \partial_{zz}^2\tilde{v} 
     =
     \int_{\Omega} (f\vec{k}\times \partial_{z}\tilde{v})\cdot \partial_{z}\tilde{v} 
     =
     0.
\end{align*}
Then, we estimate the four terms on the right side of the above equation. 
First, by H\"older's inequality and Lemma~\ref{L3}, we have 
\begin{align*}
     &\int_{\Omega}(\tilde{v}\cdot\nabla_{2})\tilde{v}\partial_{zz}^2\tilde{v}
     =
     -\int_{\Omega}(\partial_{z}\tilde{v}\cdot\nabla_{2})\tilde{v}\partial_{z}\tilde{v}
     -
     \int_{\Omega}(\tilde{v}\cdot\nabla_{2})\partial_{z}\tilde{v}\partial_{z}\tilde{v}
     \\&
     \leq
     2\int_{\Omega}|\tilde{v}| |\partial_{z}\tilde{v}| |\nabla_2\partial_{z}\tilde{v}|
     +
     \int_{\Omega}|\tilde{v}| |\partial_{z}\tilde{v}| |\partial_{zz}^2\tilde{v}|
     \\&
     \leq
     C\Vert\tilde{v}\Vert_{L^6} \Vert\partial_{z}\tilde{v}\Vert_{L^3} \Vert\nabla_2\partial_{z}\tilde{v}\Vert_{L^2}
     +
     C\Vert\tilde{v}\Vert_{L^6} \Vert\partial_{z}\tilde{v}\Vert_{L^3} \Vert\partial_{zz}^2\tilde{v}\Vert_{L^2}
     \\&
     \leq
     C\Vert\tilde{v}\Vert_{L^6} \Vert\partial_{z}\tilde{v}\Vert_{L^2}^{1/2} \Vert\nabla\partial_{z}\tilde{v}\Vert_{L^2}^{3/2}
     \\&
     \leq
     C\left(\frac{1}{\nu_{v}^3}+\frac{1}{\nu_{w}^3}\right)\Vert\tilde{v}\Vert_{L^6}^4 \Vert\partial_{z}\tilde{v}\Vert_{L^2}^{2}
     +
     \frac{\nu_{v}}{8} \Vert\nabla_{2}\partial_{z}\tilde{v}\Vert_{L^2}^{2}
     +
     \frac{\nu_{w}}{8} \Vert\partial_{zz}^2\tilde{v}\Vert_{L^2}^{2},
\end{align*}
where we integrated by parts in the first equality. 
Next, in order to bound the second term, we use the divergence-free condition and proceed as 
\begin{align*}
     &\int_{\Omega}\tilde{w}\partial_{z}\tilde{v}\partial_{zz}^2\tilde{v}
     =
     -\int_{\Omega}\left( \int_{-h}^{z}\nabla_2\cdot\tilde{v}\,d\zeta\right)\partial_{z}\tilde{v}\partial_{zz}^2\tilde{v}
     \\&
     =
     \int_{\Omega}\nabla_2\cdot\tilde{v}(\partial_{z}\tilde{v})^2
     +
     \int_{\Omega}\tilde{w}\partial_{zz}^2\tilde{v}\partial_{z}\tilde{v}.
\end{align*}
Thus, we get 
\begin{align*}
     &\int_{\Omega}\tilde{w}\partial_{z}\tilde{v}\partial_{zz}^2\tilde{v}
     =
     -\frac{1}{2}\int_{\Omega}\nabla_2\cdot\tilde{v}(\partial_{z}\tilde{v})^2
     =
     \int_{\Omega}\tilde{v}\partial_{z}\tilde{v}\nabla_2\partial_{z}\tilde{v}
     \\&
     \leq
     C\Vert\tilde{v}\Vert_{L^6} \Vert\partial_{z}\tilde{v}\Vert_{L^3} \Vert\nabla_2\partial_{z}\tilde{v}\Vert_{L^2}
     \\&
     \leq
     C\left(\frac{1}{\nu_{v}^3}+\frac{1}{\nu_{v}^2\nu_{w}}\right)\Vert\tilde{v}\Vert_{L^6}^4 \Vert\partial_{z}\tilde{v}\Vert_{L^2}^{2}
     +
     \frac{\nu_{v}}{8} \Vert\nabla_{2}\partial_{z}\tilde{v}\Vert_{L^2}^{2}
     +
     \frac{\nu_{w}}{8} \Vert\partial_{zz}^2\tilde{v}\Vert_{L^2}^{2}.
\end{align*}
For the third term, integration by parts yields 
\begin{align*}
     &-\int_{\Omega}(\nabla_{2}\int_{-h}^{z}\tilde{\theta}\,d\zeta)\cdot \partial_{zz}^2\tilde{v}
     =
     \int_{\Omega}(\nabla_{2}\tilde{\theta})\cdot \partial_{z}\tilde{v}
     =
     \int_{\Omega}\tilde{\theta}\nabla_{2}\partial_{z}\tilde{v}
     \\&
     \leq
     \Vert\tilde{\theta}\Vert_{L^2} \Vert\nabla_{2}\partial_{z}\tilde{v}\Vert_{L^2}
     \leq
     \frac{C}{\nu_{v}}\Vert\tilde{\theta}\Vert_{L^2}^2
     +
     \frac{\nu_{v}}{8} \Vert\nabla_{2}\partial_{z}\tilde{v}\Vert_{L^2}^2,
\end{align*}
where we used H\"older's and Cauchy-Schwarz inequalities. 
Regarding the last term, by condition \eqref{C2} of $I_{\delta}$ and Poincar\'e's inequality, we have 
\begin{align*}
     &-\mu_{u}\int_{\Omega}(I_{\delta}(v) - I_{\delta}(\tilde{v}))\partial_{zz}^2\tilde{v}
     \leq
     \mu_{u}\left(\Vert I_{\delta}(v)\Vert_{L^2} + \Vert I_{\delta}(\tilde{v})\Vert_{L^2}\right) \Vert\partial_{zz}^2\tilde{v}\Vert_{L^2}
     \\&
     \leq
     \frac{C\mu_{u}^2}{\nu_{w}}\left(\Vert v\Vert_{L^2}^2 + \Vert\tilde{v}\Vert_{L^2}^2\right)
     +
     \frac{\nu_{w}}{4} \Vert\partial_{zz}^2\tilde{v}\Vert_{L^2}^2. 
\end{align*}
Hence, by combining all the above estimates and using Gr\"onwall's inequality, we obtain 
$$\Vert\partial_{z}\tilde{v}(T)\Vert_{L^2} + \int_{t=0}^{t=T}\left( \nu_{v}\Vert\nabla_{2}\partial_{z}\tilde{v}\Vert_{L^2}^2 + \nu_{w}\Vert\partial_{zz}^2\tilde{v}\Vert_{L^2}^2\right)\leq K_5, $$
where $K_5$ depends on $\delta, \mu_{u}$ and $K_1$ through $K_4$.

\subsubsection{$\Vert\tilde{v}\Vert_{H^1}$ estimates}
Multiply the $\tilde{v}$-equation in \eqref{PE_DA} by $-\Delta_2\tilde{v}$, integrate over $\Omega$, and we have 
\begin{align*}
     &\frac{1}{2}\frac{d}{dt}\Vert\nabla_{2}\tilde{v}\Vert_{L^2}^2
     +
     \nu_{v}\Vert\Delta_{2}\tilde{v}\Vert_{L^2}^2
     +
     \nu_{w}\Vert\nabla_{2}\partial_{z}\tilde{v}\Vert_{L^2}^2
     \\&
     =
     \int_{\Omega}(\tilde{v}\cdot\nabla_{2})\tilde{v}\Delta_{2}\tilde{v}
     +
     \int_{\Omega}\tilde{w}\partial_{z}\tilde{v}\Delta_{2}\tilde{v}
     -
     \int_{\Omega}(\nabla_{2}\int_{-h}^{z}\tilde{\theta}\,d\zeta)\cdot \Delta_{2}\tilde{v}
     \\&\quad
     +
     \int_{\Omega}f\vec{k}\times\tilde{v}\cdot\Delta_{2}\tilde{v}
     -
     \mu_{u}\int_{\Omega}(I_{\delta}(v) - I_{\delta}(\tilde{v}))\Delta_{2}\tilde{v},
\end{align*}
where we used the equation $\partial_{z}\tilde{p} = -\tilde{\theta}$. 
We notice that integration by parts yields  
$$\int_{\Omega}f\vec{k}\times\tilde{v}\cdot\Delta_{2}\tilde{v} = 0. $$
We bound the first term on the right side of the above equation by 
\begin{align*}
     &\int_{\Omega}|\tilde{v}| |\nabla_2\tilde{v}| |\Delta_2\tilde{v}|
     \leq
     \Vert\tilde{v}\Vert_{L^6} \Vert\nabla_2\tilde{v}\Vert_{L^3} \Vert\Delta_{2}\tilde{v}\Vert_{L^2}
     \\&
     \leq
     C\Vert\tilde{v}\Vert_{L^6} \Vert\nabla_2\tilde{v}\Vert_{L^2}^{1/2} \Vert\Delta_{2}\tilde{v}\Vert_{L^2}^{3/2}
     \leq
     \frac{C}{\nu_{v}^3}\Vert\tilde{v}\Vert_{L^6}^4 \Vert\nabla_2\tilde{v}\Vert_{L^2}^2
     +
     \frac{\nu_{v}}{8}\Vert\Delta_{2}\tilde{v}\Vert_{L^2}^2.
\end{align*}
Using divergence-free condition, we estimate the second term as 
\begin{align*}
     &\int_{\Omega}\tilde{w}\partial_{z}\tilde{v}\Delta_{2}\tilde{v}
     =
     -\int_{\Omega}\left( \int_{-h}^{z}\nabla_2\cdot\tilde{v}\,d\zeta \right)\partial_{z}\tilde{v}\Delta_2\tilde{v}
     \\&
     \leq
     \int_{\mathbb{T}^2}\left( \int_{-h}^{0} |\nabla_2\tilde{v}|\,dz \int_{-h}^{0} |\partial_{z}\tilde{v}| |\Delta_2\tilde{v}|\,dz\right)\,dx\,dy
     \\&
     \leq
     C\Vert \int_{-h}^{0} |\nabla_2\tilde{v}|\,dz\Vert_{L^4} \Vert \partial_{z}\tilde{v} \Vert_{L^4} \Vert\Delta_2\tilde{v}\Vert_{L^2}
     \\&
     \leq
     \Vert\nabla_{2}\tilde{v}\Vert_{L^2}^{1/2} \Vert\partial_{z}\tilde{v}\Vert_{L^2}^{1/2} \Vert\nabla_{2}\partial_{z}\tilde{v}\Vert_{L^2}^{1/2} \Vert\Delta_{2}\tilde{v}\Vert_{L^2}^{3/2}
     \\&
     \leq
     \frac{C}{\nu_{v}^3}\left( \Vert\partial_{z}\tilde{v}\Vert_{L^2}^2\Vert\nabla_{2}\partial_{z}\tilde{v}\Vert_{L^2}^2 \right)\Vert\nabla_{2}\tilde{v}\Vert_{L^2}^2
     +
     \frac{\nu_{v}}{8}\Vert\Delta_{2}\tilde{v}\Vert_{L^2}^2, 
\end{align*}
where we used anisotropic Sobolev inequality. 
By Cauchy-Schwarz inequality, we bound the third term by 
\begin{align*}
     &\int_{\Omega} |\nabla_{2}\tilde{\theta}| |\Delta_{2}\tilde{v}|
     \leq
     \frac{C}{\nu_{v}}\Vert\nabla_{2}\tilde{\theta}\Vert_{L^2}^2
     +
     \frac{\nu_{v}}{8}\Vert\Delta_{2}\tilde{v}\Vert_{L^2}^2. 
\end{align*}
As for the last term, by condition \eqref{C2}, we have 
\begin{align*}
     &-\mu_{u}\int_{\Omega}(I_{\delta}(v) - I_{\delta}(\tilde{v}))\Delta_{2}\tilde{v}
     \\&
     \leq
     \mu_{u}\left(\Vert I_{\delta}(v)\Vert_{L^2} + \Vert I_{\delta}(\tilde{v}))\Vert_{L^2}\right) \Vert\Delta_{2}\tilde{v}\Vert_{L^2}
     \\&
     \leq
     \frac{C\mu_{u}^2}{\nu_{v}}\left(\Vert I_{\delta}(v)\Vert_{L^2}^2 + \Vert I_{\delta}(\tilde{v}))\Vert_{L^2}^2\right)
     +
     \frac{\nu_{v}}{8} \Vert\Delta_{2}\tilde{v}\Vert_{L^2}^2,
\end{align*}
where we used H\"older's and Minkowski inequalities. 
Combining all the above estimates, and by Gr\"onwall's inequality, we obtain
$$\Vert\nabla_2\tilde{v}(T)\Vert_{L^2}^2 < K_6 \quad\text{as well as}\quad\int_{0}^{T}\Vert\Delta_2\tilde{v}\Vert_{L^2}^2\,dt < \infty$$
for arbitrary $T>0$, where $K_6$ depends on $\delta, \mu_{u}$ and $K_1$ through $K_5$.. 

\subsubsection{$\Vert\tilde{\theta}\Vert_{H^1}$ estimates}
Multiplying the $\tilde{\theta}$-equation in \eqref{PE_DA} by $-(\Delta_2 + \partial_{zz}^2)\tilde{\theta}$, 
integrate over $\Omega$, 
we get 
\begin{align*}
     &\frac{1}{2}\frac{d}{dt}\left( \Vert\nabla_{2}\tilde{\theta}\Vert_{L^2}^2 + \Vert\partial_{z}\tilde{\theta}\Vert_{L^2}^2 + \Vert\nabla_{2}\tilde{\theta}(z=0)\Vert_{L^2}^2\right) 
     +
     \eta_{v}\Vert\Delta_{2}\tilde{\theta}\Vert_{L^2}^2
     +
     \eta_{w}\Vert\partial_{zz}^2\tilde{\theta}\Vert_{L^2}^2
     \\&\qquad
     +
     (\eta_{v}+\eta_{w}) (\Vert\nabla_{2}\partial_{z}\tilde{\theta}\Vert_{L^2}^2 + \Vert\nabla_{2}\tilde{\theta}(z=0)\Vert_{L^2}^2)
     \\&
     =
     \int_{\Omega}(\tilde{v}\cdot\nabla_{2})\tilde{\theta}(\Delta_2\tilde{\theta} + \partial_{zz}^2\tilde{\theta})
     -
     \int_{\Omega}\left(\int_{-h}^{z}\nabla_{2}\cdot\tilde{v}\,d\zeta\right)\partial_{z}\tilde{\theta}(\Delta_2\tilde{\theta} + \partial_{zz}^2\tilde{\theta})
     \\&\quad
     -
     \int_{\Omega}q(\Delta_2\tilde{\theta} + \partial_{zz}^2\tilde{\theta})
     -
     \mu_{\theta} \int_{\Omega}(I_{\delta}(\theta) - I_{\delta}(\tilde{\theta}))(\Delta_2\tilde{\theta} + \partial_{zz}^2\tilde{\theta}),
\end{align*}
where we used $\partial_{z}\tilde{w} = -\nabla_{2}\cdot\tilde{v}$. 
Using Lemma~\ref{L3}, the first term is bounded by 
\begin{align*}
     &\int_{\Omega}|\tilde{v}| |\nabla_{2}\tilde{\theta}| |\Delta_2\tilde{\theta} + \partial_{zz}^2\tilde{\theta}|
     \leq
     C\Vert\tilde{v}\Vert_{L^6} \Vert\nabla_{2}\tilde{\theta}\Vert_{L^3} \left(\Vert\Delta_2\tilde{\theta}\Vert_{L^2} + \Vert\partial_{zz}^2\tilde{\theta}\Vert_{L^2}\right)
     \\&
     \leq
     C\Vert\tilde{v}\Vert_{L^6} \Vert\nabla_{2}\tilde{\theta}\Vert_{L^2}^{1/2} (\Vert\Delta_{2}\tilde{\theta}\Vert_{L^2} + \Vert\nabla_{2}\partial_{z}\tilde{\theta}\Vert_{L^2})^{1/2} \left(\Vert\Delta_2\tilde{\theta}\Vert_{L^2} + \Vert\partial_{zz}^2\tilde{\theta}\Vert_{L^2}\right)
     \\&
     \leq 
     C\Vert\tilde{v}\Vert_{L^6} \Vert\nabla_{2}\tilde{\theta}\Vert_{L^2}^{1/2} \left(\Vert\Delta_{2}\tilde{\theta}\Vert_{L^2}^2 + \Vert\nabla_{2}\partial_{z}\tilde{\theta}\Vert_{L^2}^2 + \Vert\partial_{zz}^2\tilde{\theta}\Vert_{L^2}^2\right)^{3/4}
     \\&
     \leq
     C\left(\frac{1}{\eta_{v}^3}+\frac{1}{\eta_{v}^2(\eta_{v}+\eta_{w})}+\frac{1}{\eta_{v}\eta_{w}^2}+\frac{1}{\eta_{w}^2(\eta_{v}+\eta_{w})}\right)\Vert\tilde{v}\Vert_{L^6}^4 \Vert\nabla_{2}\tilde{\theta}\Vert_{L^2}^{2} 
     \\&\quad
     +
     \frac{\eta_{v}}{8}\Vert\Delta_{2}\tilde{\theta}\Vert_{L^2}^2 
     + 
     \frac{\eta_{v}+\eta_{w}}{4}\Vert\nabla_{2}\partial_{z}\tilde{\theta}\Vert_{L^2}^2 
     + 
     \frac{\eta_{w}}{8}\Vert\partial_{zz}^2\tilde{\theta}\Vert_{L^2}^2. 
\end{align*}
In order to estimate the second term, we apply the anisotropic estimates, i.e., Sobolev inequality in two-dimension, and bound it by 
\begin{align*}
     &\int_{\Omega}\left( \int_{-h}^0 |\nabla_2\tilde{v}|\,dz \int_{-h}^0 |\partial_{z}\tilde{\theta}| |\Delta_{2}\tilde{\theta} + \partial_{zz}^2\tilde{\theta}| \right)\,dx\,dy
     \\&
     \leq
     C\Vert\nabla_2\tilde{v}\Vert_{L^4} \Vert\partial_{z}\tilde{\theta}\Vert_{L^4} \Vert|\Delta_2\tilde{\theta}| + |\partial_{zz}^2\tilde{\theta}|\Vert_{L^2}
     \\&
     \leq
     C\Vert\nabla_2\tilde{v}\Vert_{L^2}^{1/2} \Vert\Delta_2\tilde{v}\Vert_{L^2}^{1/2} \Vert\partial_{z}\tilde{\theta}\Vert_{L^2}^{1/2} \left( \Vert\nabla_2\partial_{z}\tilde{\theta}\Vert_{L^2}^{2} + \Vert\partial_{zz}^2\tilde{\theta}\Vert_{L^2}^{2} \right)^{1/4} 
     \\&\quad
     \times\left( \Vert\Delta_2\tilde{\theta}\Vert_{L^2}^2 + \Vert\partial_{zz}^2\tilde{\theta}\Vert_{L^2}^2 \right)^{1/2}
     \\&
     \leq
     C\Vert\nabla_2\tilde{v}\Vert_{L^2}^{1/2} \Vert\Delta_2\tilde{v}\Vert_{L^2}^{1/2} \Vert\partial_{z}\tilde{\theta}\Vert_{L^2}^{1/2} \left( \Vert\Delta_{2}\tilde{\theta}\Vert_{L^2}^2 + \Vert\nabla_{2}\partial_{z}\tilde{\theta}\Vert_{L^2}^2 + \Vert\partial_{zz}^2\tilde{\theta}\Vert_{L^2}^2\right)^{3/4}
     \\&
     \leq
     C\left(\frac{1}{\eta_{v}^2(\eta_{v}+\eta_{w})}+\frac{1}{\eta_{v}^2\eta_{w}}+\frac{1}{\eta_{w}^2(\eta_{v}+\eta_{w})}+\frac{1}{\eta_{w}^3}\right)\Vert\nabla_2\tilde{v}\Vert_{L^2}^{2} \Vert\Delta_2\tilde{v}\Vert_{L^2}^{2} \Vert\partial_{z}\tilde{\theta}\Vert_{L^2}^{2}
     \\&\quad
     +
     \frac{\eta_{v}}{8}\Vert\Delta_{2}\tilde{\theta}\Vert_{L^2}^2 
     + 
     \frac{\eta_{v}+\eta_{w}}{4}\Vert\nabla_{2}\partial_{z}\tilde{\theta}\Vert_{L^2}^2 
     + 
     \frac{\eta_{w}}{8}\Vert\partial_{zz}^2\tilde{\theta}\Vert_{L^2}^2. 
\end{align*}
By Cauchy-Schwarz inequality, the third term on the right side of the above equation is bounded by 
\begin{align*}
     &\Vert q\Vert_{L^2} (\Vert\Delta_2\tilde{\theta}\Vert_{L^2} + \Vert\partial_{zz}^2\tilde{\theta}\Vert_{L^2})
     \\&
     \leq
     C\left(\frac{1}{\eta_{v}}+\frac{1}{\eta_{w}}\right)\Vert q\Vert_{L^2}^2
     +
     \frac{\eta_{v}}{8} \Vert\Delta_2\tilde{\theta}\Vert_{L^2}^2
     +
     \frac{\eta_{w}}{8} \Vert\partial_{zz}^2\tilde{\theta}\Vert_{L^2}^2.
\end{align*}
Regarding the last term, we use condition \eqref{C2} and H\"older's and Minkowski inequalities and get 
\begin{align*}
     &
     -\mu_{\theta} \int_{\Omega}(I_{\delta}(\theta) - I_{\delta}(\tilde{\theta}))(\Delta_2\tilde{\theta} + \partial_{zz}^2\tilde{\theta})
     \\&
     \leq
     C\mu_{\theta}\left(\Vert I_{\delta}(v)\Vert_{L^2} + \Vert I_{\delta}(\tilde{v})\Vert_{L^2}\right)(\Vert\Delta_2\tilde{\theta}\Vert_{L^2} + \Vert\partial_{zz}^2\tilde{\theta}\Vert_{L^2})
     \\&
     \leq
     \frac{C\mu_{\theta}^2}{\eta_{v}}\left(\Vert I_{\delta}(v)\Vert_{L^2}^2 + \Vert I_{\delta}(\tilde{v})\Vert_{L^2}^2\right)
     +
     \frac{C\mu_{\theta}^2}{\eta_{w}}\left(\Vert I_{\delta}(v)\Vert_{L^2}^2 + \Vert I_{\delta}(\tilde{v})\Vert_{L^2}^2\right)
     \\&\quad
     +
     \frac{\eta_{v}}{8}\Vert\Delta_2\tilde{\theta}\Vert_{L^2}^2
     +
     \frac{\eta_{w}}{8}\Vert\partial_{zz}^2\tilde{\theta}\Vert_{L^2}^2. 
\end{align*}
Summing up all the above estimates, and using the $H^1$ estimates for $\tilde{v}$, as well as Gr\"onwall's inequality, we get
$$\Vert\nabla\tilde{\theta}(T)\Vert_{L^2}^2 < K_7$$
for arbitrary $T>0$, where $K_7$ depends on $\Vert q\Vert_{H^1}$ and $K_1$ through $K_6$. 

Now, with the above $L^2$, $L^6$, and $H^1$ estimates ready, 
suppose that there exists some time $T_{\ast}>0$ such that the interval $(0, T_{\ast})$ is the maximal interval of existence of the solution $(\tilde{u}, \tilde{\theta}, \tilde{p})$. 
We infer that 
$$\lim_{t \to T_{\ast}^{-}}\sup{\{\Vert\tilde{v}(t)\Vert_{H^1} + \Vert\tilde{\theta}(t)\Vert_{H^1}\}} = \infty$$
must be valid provided $T_{\ast} < \infty$. 
Therefore, we obtain a contradiction to the above estimates and the solution thus exists globally for all $T>0$. 
Proof of Theorem~\ref{T0} is complete. 
\end{proof}

\section{Exponential-rate Convergence to the Reference Solution in $L^2$-norm}\label{Exp_Decay}
In this section, we prove Theorem~\ref{T4}. 
\begin{proof}

By subtracting the corresponding equations of system \eqref{PE_DA} from those of system \eqref{PE}, we obtain 
\begin{equation*}
 \left\{
   \begin{aligned}
   &
   \frac{\partial V}{\partial t} + (V\cdot\nabla_{2})v + W\frac{\partial v}{\partial z} 
   +
   ((v-V)\cdot\nabla_2)V + (w-W)\frac{\partial V}{\partial z} 
   \\&
   =
   -
   \nabla_{2} P 
   - 
   f\vec{k}\times V 
   -L_{1}V 
   -
   \mu_{u}I_{\delta}(V), 
   \\&
   \partial_{z}P + \Theta = 0,
   \\&
   \nabla_{2}\cdot V + \partial_{z}W = 0,
   \\&
   \frac{\partial \Theta}{\partial t} + V\cdot\nabla_{2} \theta + W\frac{\partial \theta}{\partial z} 
   + 
   (v-V)\cdot\nabla_{2}\Theta + (w-W)\frac{\partial\Theta}{\partial z}
   \\&
   = 
   -L_{2}\theta 
   - 
   \mu_{\theta}I_{\delta}\Theta,
   \\&
   \text{with} \quad V(x,0) = v(x,0)-0 = v_0(x), \quad \text{and} \quad \Theta(x,0) = \theta(x, 0)-0 = \theta_{0}(x), 
  \end{aligned}
 \right.
\end{equation*}
where we denote by $U=(V, W)$ the difference $V=v-\tilde{v}$, $W=w-\tilde{w}$, $\Theta=\theta-\tilde{\theta}$, and $P=p-\tilde{p}$. 
For $-h<z<0$, integrating the equation of the pressure $P$ vertically, we get 
$$
     P = -\int_{-h}^{z} \Theta\,d\zeta + P_0, 
$$
and without loss of generality, we assume $P_0 = 0$. 
Then, we multiply the equations of $V$ and $\Theta$ by $V$ and $\Theta$, respectively, use the above equation of $P$, integrate over $\Omega$, and get 
\begin{equation*}
 \left\{
   \begin{aligned}
   &
   \frac{1}{2}\frac{d}{d t}\Vert V\Vert_{L^2}^2 
   +
   \nu_{v}\Vert\nabla_{2}V\Vert_{L^2}^2 + \nu_{w}\Vert\partial_{z}V\Vert_{L^2}^2
   \\&
   = 
   -\int_{\Omega}(V\cdot\nabla_{2})v\cdot V 
   +
   \int_{\Omega}\left(\int_{-h}^{z}\nabla_{2}\cdot V\,d\zeta\right)\partial_{z}v\cdot V 
   \\&\qquad
   -
   \int_{\Omega}\left((v-V)\cdot\nabla_2\right)V\cdot V 
   +
   \int_{\Omega}\left(\int_{-h}^{z}-\nabla_{2}(v-V)\,d\zeta\right)\partial_{z}V\cdot V 
   \\&\qquad
   + 
   \int_{\Omega}(\nabla_{2}\int_{-h}^{z}\Theta\,d\zeta)\cdot V 
   -
   \mu_{u}\int_{\Omega}I_{\delta}(V)\cdot V, 
   \\&
   \frac{1}{2}\frac{d}{d t}\Vert \Theta\Vert_{L^2}^2 
   +
   \eta_{v}\Vert\nabla_{2}\Theta\Vert_{L^2}^2 + \eta_{w}\Vert\partial_{z}\Theta\Vert_{L^2}^2 
   \\&
   = 
   -\int_{\Omega}\Theta\left(V\cdot\nabla_{2} \theta\right)
   + 
   \int_{\Omega}\left(\int_{-h}^{z}\nabla_{2}\cdot V\,d\zeta\right)\Theta\partial_{z}\theta 
   \\&\qquad
   - 
   \int_{\Omega}\Theta\left( (v-V)\cdot\nabla_{2}\Theta\right) 
   + 
   \int_{\Omega}\Theta\left(\int_{-h}^{z}-\nabla_{2}(v-V)\,d\zeta\right)\partial_{z}\Theta
   \\&\qquad
   - 
   \mu_{\theta}\int_{\Omega}\Theta I_{\delta}(\Theta),
  \end{aligned}
 \right.
\end{equation*}
where we used the divergence-free condition $\partial_{z}W = -\nabla_{2}\cdot V$. 
Next, we estimate the right sides of the above two equation of $V$ and $\Theta$. 
Once again in view of the divergence-free condition, we deduce that the third and fourth terms on the right side of the $V$-equation add up to zero. 
For the same reason, the third and fourth terms on the right side of the $\Theta$-equation also sum up to $0$. 
In order to bound the first term on the right side of the $V$-equation, we use H\"older's inequality and Lemma~\ref{L3}, and obtain  
\begin{align*}
     &
     -\int_{\Omega}(V\cdot\nabla_{2})v\cdot V
     \leq
     \Vert\nabla_{2}v\Vert_{L^2} \Vert V\Vert_{L^3} \Vert V\Vert_{L^6}
     \leq
     \Vert\nabla_{2}v\Vert_{L^2} \Vert V\Vert_{L^2}^{1/2} \Vert\nabla_{2}V\Vert_{L^2}^{3/2}
     \\&
     \leq
     \frac{C}{\nu_{v}^3}\Vert\nabla_{2}v\Vert_{L^2}^{4} \Vert V\Vert_{L^2}^{2}
     + 
     \frac{\nu_{v}}{8}\Vert\nabla_{2}V\Vert_{L^2}^{2}.
\end{align*}
For the second term on the right side of the $V$-equation, we estimate after integration by parts as 
\begin{align*}
     &
     \int_{\Omega}\left(\int_{-h}^{z}\nabla_{2}\cdot V\,d\zeta\right)\partial_{z}v\cdot V
     \\&
     =
     -\sum_{i=1}^{2}\int_{\Omega}\left(\int_{-h}^{z}V\,d\zeta\right)\partial_{iz}^{2}v\cdot V 
     -
     \sum_{i=1}^{2}\int_{\Omega}\left(\int_{-h}^{z}V\,d\zeta\right)\partial_{z}v\cdot \partial_{i}V
     \\&
     \leq
     \Vert\Delta v\Vert_{L^2} \Vert V\Vert_{L^3} \Vert V\Vert_{L^6}
     +
     \Vert\partial_{z}v\Vert_{L^6} \Vert V\Vert_{L^3} \Vert\nabla V\Vert_{L^2}
     \\&
     \leq
     2\Vert\Delta v\Vert_{L^2} \Vert V\Vert_{L^2}^{1/2} \Vert\nabla V\Vert_{L^2}^{3/2}
     \\&
     \leq
     C\left(\frac{1}{\nu_{v}^3}+\frac{1}{\nu_{w}^3}\right)\Vert\Delta v\Vert_{L^2}^4 \Vert V\Vert_{L^2}^2 
     +
     \frac{\nu_{v}}{8}\Vert\nabla_{2}V\Vert_{L^2}^2
     +
     \frac{\nu_{w}}{8}\Vert\partial_{z}V\Vert_{L^2}^2. 
\end{align*}
Regarding the last but second term on the right side of the $V$-equation, we have 
\begin{align*}
     &
     \int_{\Omega}(\nabla_{2}\int_{-h}^{z}\Theta\,d\zeta)\cdot V
     =
     -\int_{\Omega}(\int_{-h}^{z}\Theta\,d\zeta)\nabla_{2}\cdot V
     \leq
     C\Vert\Theta\Vert_{L^2} \Vert\nabla_{2}V\Vert_{L^2}
     \\&
     \leq
     \frac{C}{\nu_{v}}\Vert\Theta\Vert_{L^2}^2 
     +
     \frac{\nu_{v}}{8}\Vert\nabla_{2}V\Vert_{L^2}^2,
\end{align*}
where we integrated by parts and used Cauchy-Schwarz inequality. 
The last term in the $V$-equation is bounded as 
\begin{align*}
     &
     -\mu_{u}\int_{\Omega}I_{\delta}(V)\cdot V
     =
     -\mu_{u}\int_{\Omega}(I_{\delta}(V) - V)\cdot V
     -
     \mu_{u}\Vert V\Vert_{L^2}^2
     \\&
     \leq
     c_1\mu_{u}\delta\Vert\nabla V\Vert_{L^2}\Vert V\Vert_{L^2}
     -
     \mu_{u}\Vert V\Vert_{L^2}^2
     \\&
     \leq
     C\left(\frac{1}{\nu_{v}}+\frac{1}{\nu_{w}}\right)c_{1}^2\delta^2\mu_{u}^2 \Vert V\Vert_{L^2}^2
     +
     \frac{\nu_{v}}{8}\Vert\nabla_{2}V\Vert_{L^2}^2
     +
     \frac{\nu_{w}}{8}\Vert\partial_{z}V\Vert_{L^2}^2
     -
     \mu_{u}\Vert V\Vert_{L^2}^2.
\end{align*}

Next, we estimate the right side of the $\Theta$-equation. 
For the first term, by H\"older's inequality and Lemma~\ref{L3}, we get 
\begin{align*}
     &
     \int_{\Omega}\Theta\left(V\cdot\nabla_{2} \theta\right)
     \leq
     \Vert\nabla_{2}\theta\Vert_{L^6} \Vert V\Vert_{L^2} \Vert \Theta\Vert_{L^3}
     \leq
     \Vert\Delta\theta\Vert_{L^2} \Vert V\Vert_{L^2} \Vert\Theta\Vert_{L^2}^{1/2} \Vert\nabla\Theta\Vert_{L^2}^{1/2}
     \\&
     \leq
     C\Vert\Delta\theta\Vert_{L^2}^2 \Vert V\Vert_{L^2}^2
     +
     C\left(\frac{1}{\nu_{v}}+\frac{1}{\nu_{w}}\right)\Vert\Theta\Vert_{L^2}^2
     \\&\quad
     +
     \frac{\eta_{v}}{8}\Vert\nabla_{2}\Theta\Vert_{L^2}^2
     +
     \frac{\eta_{w}}{8}\Vert\partial_{z}\Theta\Vert_{L^2}^2. 
\end{align*}
In order to estimate the second term, we first integrate by parts, then proceed as 
\begin{align*}
     &
     \int_{\Omega}\left(\int_{-h}^{z}\nabla_{2}\cdot V\,d\zeta\right)\Theta\partial_{z}\theta
     \\&
     =
     -\sum_{i=1}^{2}\int_{\Omega}\left(\int_{-h}^{z}V\,d\zeta\right)\Theta\partial_{iz}^2\theta
     -
     \sum_{i=1}^{2}\int_{\Omega}\left(\int_{-h}^{z}V\,d\zeta\right)\partial_{i}\Theta\partial_{z}\theta
     \\&
     \leq
     \Vert\Delta\theta\Vert_{L^2} \Vert V\Vert_{L^3} \Vert\Theta\Vert_{L^6}
     +
     \Vert\partial_{z}\theta\Vert_{L^6} \Vert V\Vert_{L^3} \Vert\nabla\Theta\Vert_{L^2}
     \\&
     \leq
     2\Vert\Delta\theta\Vert_{L^2} \Vert V\Vert_{L^2}^{1/2} \Vert\nabla V\Vert_{L^2}^{1/2} \Vert\nabla\Theta\Vert_{L^2}
     \\&
     \leq
     C\left(\frac{1}{\nu_{v}\eta_{v}^2}+\frac{1}{\nu_{w}\eta_{v}^2}+\frac{1}{\nu_{v}\eta_{w}^2}+\frac{1}{\nu_{w}\eta_{w}^2}\right)\Vert\Delta\theta\Vert_{L^2}^4 \Vert V\Vert_{L^2}^{2}
     \\&\quad
     +
     \frac{\nu_{v}}{8}\Vert\nabla_{2}V\Vert_{L^2}^2
     +
     \frac{\nu_{w}}{8}\Vert\partial_{z}V\Vert_{L^2}^2
     +
     \frac{\eta_{v}}{8}\Vert\nabla_{2}\Theta\Vert_{L^2}^2
     +
     \frac{\eta_{w}}{8}\Vert\partial_{z}\Theta\Vert_{L^2}^2.
\end{align*}
As for the last term, we bound it as 
\begin{align*}
     &
     -\mu_{\theta}\int_{\Omega}\Theta I_{\delta}(\Theta)
     =
     -\mu_{\theta}\int_{\Omega}\Theta (I_{\delta}(\Theta) - \Theta)
     -
     \mu_{\theta}\Vert\Theta\Vert_{L^2}^2
     \\&
     \leq
     c_1\mu_{\theta}\delta\Vert\Theta\Vert_{L^2}\Vert\nabla\Theta\Vert_{L^2}
     -
     \mu_{\theta}\Vert\Theta\Vert_{L^2}^2
     \\&
     \leq
     C\left(\frac{1}{\eta_{v}}+\frac{1}{\eta_{w}}\right)c_{1}^2\delta^2\mu_{\theta}^2 \Vert\Theta\Vert_{L^2}^2
     +
     \frac{\eta_{v}}{8}\Vert\nabla_{2}\Theta\Vert_{L^2}^2
     +
     \frac{\eta_{w}}{8}\Vert\partial_{z}\Theta\Vert_{L^2}^2
     -
     \mu_{\theta}\Vert\Theta\Vert_{L^2}^2. 
\end{align*}
Combining all the above estimates, and by denoting 
$$X(t) = \Vert V\Vert_{L^2}^2 + \Vert\Theta\Vert_{L^2}^2 \quad\text{ and }\quad Y(t) = \Vert\nabla V\Vert_{L^2}^2 + \Vert\nabla\Theta\Vert_{L^2}^2, $$
and using Theorem~\ref{T1} and Theorem~\ref{T2}, we obtain 
\begin{align}
     &\frac{d}{dt}X(t) + \frac{1}{2}Y(t)
     \\&
     \leq
     \left(C + C_{\nu, \eta}\bigl(\Vert\nabla_{2}v\Vert_{L^2}^4 + \Vert\Delta v\Vert_{L^2}^4 
     + \Vert\nabla\theta\Vert_{L^2}^4 + \Vert\Delta\theta\Vert_{L^2}^2 
     + \Vert\Delta\theta\Vert_{L^2}^4 \bigr) \right. \nonumber
     \\&\qquad
     \left.+ \left(\frac{1}{\nu_{v}}+\frac{1}{\nu_{w}}\right)c_1^2\delta^2\mu_{u}^2
     + \left(\frac{1}{\eta_{v}}+\frac{1}{\eta_{w}}\right)c_1^2\delta^2\mu_{\theta}^2
     - \mu_{min} \right)X(t), \nonumber
     \label{L2_exp}
\end{align}
where $\mu_{min} = \min\{\mu_{u}, \mu_{\theta}\}$ and C depends on the reference solution, $\lambda_1$, and the domain. 
Now choose 
$$0 < \delta < \frac{1}{\sqrt{2}c_1}\frac{\sqrt{\mu_{min}}}{\mu_{max}\sqrt{\nu_{v}^{-1}+\nu_{w}^{-1}+\eta_{v}^{-1}+\eta_{w}^{-1}}}$$
sufficiently small such that 
$$\left(\frac{1}{\nu_{v}}+\frac{1}{\nu_{w}}\right)c_1^2\delta^2\mu_{u}^2
     + \left(\frac{1}{\eta_{v}}+\frac{1}{\eta_{w}}\right)c_1^2\delta^2\mu_{\theta}^2 < \frac{\mu_{min}}{2}$$
where $\mu_{max} = \max\{\mu_{u}, \mu_{\theta}\}$, 
and choose 
$$\mu_{min} > 2\left(C + C_{\nu, \eta}(\Vert v\Vert_{H^2}^4 + \Vert\theta\Vert_{H^2}^4)\right), $$
such that the coefficient of $X(t)$ on the right side of \eqref{L2_exp} is negative, 
where the constant $C_{\nu, \eta}$ also depends on $\nu_{v}, \nu_{w}, \eta_{v}, \eta_{w}$. 
Then, by Gr\"onwall's inequality, we conclude that $X(t)$ tends to zero exponentially fast as time $t \to \infty$. 
Thus, the exponential $L^2$-convergence of the data assimilation solution to the reference solution is proven. 
\end{proof}

\begin{remark}
     In view of the available results regarding the regularity of solutions to the 3D primitive equations and the continuous data assimilation on other dynamical systems, we believe that the $H^1$ norms of $V$ and $\Theta$ should also tend to zero at an exponential rate and the relevant work will be presented in future works. 
\end{remark}


\begin{thebibliography}{99}

\bibitem{EB}
     \newblock  R. Errico and D. Baumhefner, 
     \newblock Predictability experiments using a high-resolution limited-area model, 
     \newblock \emph{Monthly Weather Review}, \textbf{115} (1986), 488--505.
     
\bibitem{HA}
     \newblock  J. Hoke and R. Anthes, 
     \newblock The initialization of numerical models by a dynamic relaxation technique, 
     \newblock \emph{Monthly Weather Review}, \textbf{104} (1976), 1551--1556.

\bibitem{LAE}
     \newblock  A. Lorenc, W. Adams, and J. Eyre, 
     \newblock The treatment of humidity in ECMWF's data assimilation scheme, Atmospheric Water Vapor, 
     \newblock \emph{Academic Press New York}, (1980), 497--512.
     
\bibitem{AOT} (MR3183055) [10.1007/s00332-013-9189-y]
    \newblock A. Azouani, E. Olson, and E. S. Titi, 
    \newblock Continuous data assimilation using general interpolant observables, 
    \newblock \emph{J. Nonlinear Sci.}, \textbf{24} (2014), 277--304.

\bibitem{AT} (MR3274649) [10.3934/eect.2014.3.579]
    \newblock A. Azouani and E. S. Titi, 
    \newblock Feedback control of nonlinear dissipative systems by finite determining
	parameters---a reaction-diffusion paradigm, 
    \newblock \emph{Evol. Equ. Control Theory}, \textbf{3} (2014), 579--594.

\bibitem{CKT} (MR1855664) [10.1512/iumj.2001.50.2154]
    \newblock C. Cao, I. G. Kevrekidis, and E. S. Titi, 
    \newblock Numerical criterion for the stabilization of steady states of the {N}avier--{S}tokes equations,  
    \newblock \emph{Indiana Univ. Math. J.}, \textbf{50} (2001), 37--96.

\bibitem{HOT} (MR2831793) [10.1016/j.physd.2011.04.021]
    \newblock K. Hayden, E. Olson, and E. S. Titi, 
    \newblock Discrete data assimilation in the {L}orenz and 2{D} {N}avier-{S}tokes equations,
    \newblock \emph{Phys. D}, \textbf{240} (2011), 1416--1425.

\bibitem{OT} (MR2036872 (2004k:76059)) [10.1023/A:1027312703252]
    \newblock E. Olson and E. S. Titi, 
    \newblock Determining modes for continuous data assimilation in 2{D} turbulence,
    \newblock \emph{J. Statist. Phys.}, \textbf{113} (2003), 799--840.

\bibitem{BHLP}
    \newblock A. Biswas, J. Hudson, A. Larios, and Y. Pei, 
    \newblock Continuous data assimilation for the magneto-hydrodynamic equations in 2{D} using one component of the velocity and magnetic fields,
    \newblock \emph{Asymptot. Anal.}, \textbf{to appear}, (2018). 
    
\bibitem{ANLT} (MR3475121) [10.3233/ASY-151351]
    \newblock D. A. Albanez, H. J. Nussenzveig Lopes, E. S. Titi, 
    \newblock Continuous data assimilation for the three-dimensional {N}avier--{S}tokes-$\alpha$ model, 
    \newblock \emph{Asymptotic Anal.}, \textbf{97} (2016), 139--164.  
    
\bibitem{ATG} (MR3650450) [10.1007/s10596-017-9619-2]
    \newblock M. U. Altaf, E. S. Titi, T. Gebrael, O. M. Knio, L. Zhao, M. F. McCabe, I. Hoteit, 
    \newblock Downscaling the 2{D} {B\'e}nard convection equations using continuous data assimilation, 
    \newblock \emph{Computational Geosciences}, (2017), 1--18.      
    
\bibitem{BOT} (MR3319381) [10.1088/0951-7715/28/3/729]
    \newblock H. Bessaih, E. Olson, and E. S. Titi, 
    \newblock Continuous data assimilation with stochastically noisy data,
    \newblock \emph{Nonlinearity}, \textbf{28} (2015), 729--753. 
    
\bibitem{BM} (MR3595320) [10.1016/j.nonrwa.2016.10.005]
    \newblock A. Biswas and V. R. Martinez, 
    \newblock Higher-order synchronization for a data assimilation algorithm for the 2{D} {N}avier--{S}tokes equations,
    \newblock \emph{Nonlinearity}, \textbf{35} (2017), 132--157.    

\bibitem{FJT} (MR3349518) [10.1016/j.physd.2015.03.011]
    \newblock A. Farhat, M. Jolly, and E. S. Titi, 
    \newblock Continuous data assimilation for the 2{D} {B}\'enard convection through velocity measurements alone,
    \newblock \emph{Phys. D}, \textbf{303} (2015), 59--66. 

\bibitem{FLT1} (MR3461924) [10.1007/s00021-015-0225-6]
    \newblock A. Farhat, E. Lunasin, and E. S. Titi, 
    \newblock Abridged continuous data assimilation for the 2{D} {N}avier--{S}tokes Equations utilizing measurements of only one component of the velocity field,
    \newblock \emph{J. Math. Fluid Mech.}, \textbf{18} (2016), 1--23. 

\bibitem{FLT2} (MR3462588) [10.1016/j.jmaa.2016.01.072]
    \newblock A. Farhat, E. Lunasin, and E. S. Titi, 
    \newblock Data assimilation algorithm for 3{D} {B\'e}nard convection in porous
	media employing only temperature measurements,
    \newblock \emph{J. of Math. Anal. and Appl.}, \textbf{438} (2016), 492--506. 
    
\bibitem{FLT3} (MR3638329) [10.1007/s00332-017-9360-y]
    \newblock A. Farhat, E. Lunasin, and E. S. Titi, 
    \newblock Continuous Data Assimilation for a {2D} {B}\'enard Convection System
	Through Horizontal Velocity Measurements Alone,
    \newblock \emph{J. of Nonlinear Science}, (2017), 1--23. 
    
\bibitem{FMT} (MR3570275) [10.1137/16M1076526]
    \newblock C. Foias, C. Mondaini, and E. S. Titi, 
    \newblock A discrete data assimilation scheme for the solutions of the two-dimensional
	{N}avier-{S}tokes equations and their statistics,
    \newblock \emph{SIAM J. Appl. Dyn. Syst.}, \textbf{15} (2016), 2109--2142. 
    
\bibitem{GOT} (MR3514248) 
    \newblock M. Gesho, E. Olson, and E. S. Titi, 
    \newblock A computational study of a data assimilation algorithm for the two-dimensional {N}avier-{S}tokes equations,
    \newblock \emph{Commun. Comput. Phys.}, \textbf{19} (2016), 1094--1110.     
    
\bibitem{JMT} (MR3604950) [10.1515/ans-2016-6019]
    \newblock M. Jolly, and V. R. Martinez, and E. S. Titi,  
    \newblock A data assimilation algorithm for the subcritical surface quasi-geostrophic equation,
    \newblock \emph{Adv. Nonlinear Stud.}, \textbf{17} (2017), 167--192.    
    
\bibitem{MTT} (MR3476509) [10.1088/0951-7715/29/4/1292]
    \newblock P. A. Markowich, E. S. Titi, and S. Trabelsi,  
    \newblock Continuous data assimilation for the three-dimensional {B}rinkman--{F}orchheimer-extended {D}arcy model,
    \newblock \emph{Nonlinearity}, \textbf{29} (2016), 1292--1328.    
    
\bibitem{LTW1} (MR1158375) [10.1088/0951-7715/29/4/1292]
    \newblock J.-L. Lions, R. Temam, and S. Wang,  
    \newblock Continuous data assimilation for the three-dimensional {B}rinkman--{F}orchheimer-extended {D}arcy model,
    \newblock \emph{Nonlinearity}, \textbf{5} (1992), 237--288. 
    
\bibitem{LTW2} (MR1187737 (93k:86004)) 
    \newblock J.-L. Lions, R. Temam, and S. Wang,  
    \newblock On the equations of the large-scale ocean,
    \newblock \emph{Nonlinearity}, \textbf{5} (1992), 1007--1053. 
    
\bibitem{LTW3} (MR1325825) 
    \newblock J.-L. Lions, R. Temam, and S. Wang,  
    \newblock Mathematical theory for the coupled atmosphere-ocean models.
              ({CAO} {III}),
    \newblock \emph{J. Math. Pures Appl. (9)}, \textbf{74} (1992), 105--163.     
 
\bibitem{Z1} (MR1355151) [10.1016/0893-9659(94)00110-X]
    \newblock M. Ziane,  
    \newblock Regularity results for {S}tokes type systems related to climatology,
    \newblock \emph{Appl. Math. Lett.}, \textbf{8} (1995), 53--58. 

\bibitem{Z2} (MR1418137) [10.1016/0362-546X(95)00154-N]
    \newblock M. Ziane,  
    \newblock Regularity results for the stationary primitive equations of the atmosphere and the ocean,
    \newblock \emph{Nonlinear Anal.}, \textbf{28} (1997), 289-313. 

\bibitem{BGMR2} (MR2025211)
    \newblock D. Bresch, F. Guill\'en-Gonz\'alez, N. Masmoudi, and M. A. Rodr\'iguez-Bellido,
    \newblock  Asymptotic derivation of a {N}avier condition for the primitive equations, 
    \newblock  \emph{Asymptot. Anal.}, \textbf{33} (2003), 237-259
    
\bibitem{BGMR3} (MR1948873)
    \newblock D. Bresch, F. Guill\'en-Gonz\'alez, N. Masmoudi, and M. A. Rodr\'iguez-Bellido,
    \newblock  On the uniqueness of weak solutions of the two-dimensional primitive equations, 
    \newblock \emph{Differential Integral Equations}, \textbf{16} (2003), 77--94. 
    
\bibitem{GMR} (MR1859612)
    \newblock F. Guill\'en-Gonz\'alez, N. Masmoudi, and M. A. Rodr\'iguez-Bellido,
    \newblock  Anisotropic estimates and strong solutions of the primitive equations, 
    \newblock \emph{Differential Integral Equations}, \textbf{14} (2001), 1381--1408.  
    
\bibitem{CT1} (MR2342696) [10.4007/annals.2007.166.245]
    \newblock C. Cao and E. S. Titi,  
    \newblock Global well-posedness of the three-dimensional viscous primitive equations of large scale ocean and atmosphere dynamics,
    \newblock \emph{Ann. of Math. (2)}, \textbf{166} (2007), 245--267.
    
\bibitem{CT2} (MR2890308) [10.1007/s00220-011-1409-4]
    \newblock C. Cao and E. S. Titi,  
    \newblock Global well-posedness of the {$3D$} primitive equations with partial vertical turbulence mixing heat diffusion,
    \newblock \emph{Comm. Math. Phys.}, \textbf{310} (2012), 537--563.
    
\bibitem{CLT1} (MR3237881) [10.1007/s00205-014-0752-y]
    \newblock C. Cao, J. Li, and E. S. Titi,  
    \newblock Local and global well-posedness of strong solutions to the 3{D} primitive equations with vertical eddy diffusivity,
    \newblock \emph{Arch. Ration. Mech. Anal.}, \textbf{214} (2014), 35--76.    
    
\bibitem{CLT2} (MR3264417) [10.1016/j.jde.2014.08.003]
    \newblock C. Cao, J. Li, and E. S. Titi,  
    \newblock Global well-posedness of strong solutions to the 3{D} primitive equations with horizontal eddy diffusivity,
    \newblock \emph{J. Differential Equations}, \textbf{257} (2014), 4108--4132.    

\bibitem{CLT3} (MR3518238) [10.1002/cpa.21576]
    \newblock C. Cao, J. Li, and E. S. Titi,  
    \newblock Global well-posedness of the three-dimensional primitive equations with only horizontal viscosity and diffusion,
    \newblock \emph{Comm. Pure Appl. Math.}, \textbf{69} (2016), 1492--1531. 
    
\bibitem{CLT4} (MR3630635) [10.1016/j.jfa.2017.01.018]
    \newblock C. Cao, J. Li, and E. S. Titi,  
    \newblock Strong solutions to the 3{D} primitive equations with only horizontal dissipation: near {$H^1$} initial data,
    \newblock \emph{J. Funct. Anal.}, \textbf{272} (2017), 4606--4641. 

\bibitem{CINT} (MR3339156) [10.1007/s00220-015-2365-1]
    \newblock C. Cao, S. Ibrahim, K. Nakanishi, and E. S. Titi,  
    \newblock Finite-time blowup for the inviscid primitive equations of oceanic and atmospheric dynamics,
    \newblock \emph{Comm. Math. Phys.}, \textbf{337} (2015), 473--482. 
    
\bibitem{GH} (MR2509995) [10.1016/S0252-9602(09)60074-6]
    \newblock B. Guo and D. Huang,  
    \newblock On the 3{D} viscous primitive equations of the large-scale atmosphere,
    \newblock \emph{Acta Math. Sci. Ser. B Engl. Ed.}, \textbf{29} (2009), 846--866. 

\bibitem{GHZ} (MR2434911) [10.3934/dcdsb.2008.10.801]
    \newblock N. Glatt-Holtz and M. Ziane,  
    \newblock The stochastic primitive equations in two space dimensions with multiplicative noise,
    \newblock \emph{Discrete Contin. Dyn. Syst. Ser. B}, \textbf{10} (2008), 801--822. 

\bibitem{H} (MR2123085) [10.1016/j.na.2004.12.005]
    \newblock C. Hu,  
    \newblock Asymptotic analysis of the primitive equations under the small depth assumption,
    \newblock \emph{Nonlinear Anal.}, \textbf{61} (2005), 425--460. 

\bibitem{K} (MR2374160) [10.1007/s00021-006-0228-4]
    \newblock G. M. Kobelkov,  
    \newblock Existence of a solution ``in the large'' for ocean dynamics equations,
    \newblock \emph{J. Math. Fluid Mech.}, \textbf{9} (2007), 588--610. 

\bibitem{LT} (MR3592073) [10.1137/15M1050513]
    \newblock J. Li and E. S. Titi,  
    \newblock Existence and uniqueness of weak solutions to viscous primitive equations for a certain class of discontinuous initial data,
    \newblock \emph{SIAM J. Math. Anal.}, \textbf{49} (2017), 1--28. 

\bibitem{KPRZ} (MR3207927) [10.1088/0951-7715/27/6/1135]
    \newblock I. Kukavica and M. Ziane,  
    \newblock Primitive equations with continuous initial data,
    \newblock \emph{Nonlinearity}, \textbf{27} (2014), 1135--1155. 
    
\bibitem{KZ1} (MR2353676) [10.1016/j.crma.2007.07.025]
    \newblock I. Kukavica and M. Ziane,  
    \newblock The regularity of solutions of the primitive equations of the ocean in space dimension three,
    \newblock \emph{C. R. Math. Acad. Sci. Paris}, \textbf{345} (2007), 257--260. 
    
\bibitem{KZ2} (MR2368323 (2008k:35379)) [10.1088/0951-7715/20/12/001]
    \newblock I. Kukavica and M. Ziane,  
    \newblock On the regularity of the primitive equations of the ocean,
    \newblock \emph{Nonlinearity}, \textbf{20} (2007), 2739--2753. 
    
\bibitem{KZ3} (MR2483337) 
    \newblock I. Kukavica and M. Ziane,   
    \newblock Uniform gradient bounds for the primitive equations of the ocean,
    \newblock \emph{Differential Integral Equations}, \textbf{21} (2008), 837--849.     
    
\bibitem{P} (MR2082773) 
    \newblock M. Petcu,  
    \newblock Gevrey class regularity for the primitive equations in space dimension $2$,
    \newblock \emph{Asymptot. Anal.}, \textbf{39} (2004), 1--13.    
        
\bibitem{SV} (MR1689401) 
    \newblock M. Schonbek and G. K. Vallis,   
    \newblock Energy decay of solutions to the {B}oussinesq, primitive, and planetary geostrophic equations,
    \newblock \emph{J. Math. Anal. Appl.}, \textbf{234} (1999), 457--481.     

\bibitem{KTVZ1} (MR2652489) [10.1006/jmaa.1999.6354]
    \newblock I. Kukavica, R. Temam, V. Vicol, and M. Ziane,   
    \newblock Existence and uniqueness of solutions for the hydrostatic {E}uler equations on a bounded domain with analytic data,
    \newblock \emph{C. R. Math. Acad. Sci. Paris}, \textbf{348} (2010), 639--645.     
    
\bibitem{KTVZ2} (MR2737223) [10.1016/j.jde.2010.07.032]
    \newblock I. Kukavica, R. Temam, V. Vicol, and M. Ziane,   
    \newblock Local existence and uniqueness for the hydrostatic {E}uler equations on a bounded domain,
    \newblock \emph{J. of Differential Equations}, \textbf{250} (2011), 1719--1746.     
    
\bibitem{MW} (MR2898740) [10.1007/s00205-011-0485-0]
    \newblock N. Masmoudi and T. K. Wong,   
    \newblock On the {$H^s$} theory of hydrostatic {E}uler equations,
    \newblock \emph{Arch. Ration. Mech. Anal.}, \textbf{204} (2012), 231--271.

\bibitem{R} (MR2563627) [10.1007/s00205-008-0207-4]
    \newblock M. Renardy,   
    \newblock Ill-posedness of the hydrostatic {E}uler and {N}avier-{S}tokes equations,
    \newblock \emph{Arch. Ration. Mech. Anal.}, \textbf{194} (2009), 877--886.
    
\bibitem{RTT1} (MR2166669) [10.3934/dcds.2005.13.1257]
    \newblock A. Rousseau, A. R. Temam, and J. Tribbia,   
    \newblock Boundary conditions for the 2{D} linearized {PE}s of the ocean in the absence of viscosity,
    \newblock \emph{Discrete Contin. Dyn. Syst.}, \textbf{13} (2005), 1257--1276.
    
\bibitem{RTT2} (MR2401691) [10.1016/j.matpur.2007.12.001]
    \newblock A. Rousseau, A. R. Temam, and J. Tribbia,   
    \newblock The 3{D} primitive equations in the absence of viscosity: boundary conditions and well-posedness in the linearized case,
    \newblock \emph{J. Math. Pures Appl. (9)}, \textbf{89} (2008), 297--319.
    
\bibitem{C} (MR3233751) [10.1017/S0308210512001953]
    \newblock I. Chueshov,   
    \newblock A squeezing property and its applications to a description of long-time behaviour in the three-dimensional viscous primitive equations,
    \newblock \emph{Proc. Roy. Soc. Edinburgh Sect. A}, \textbf{144} (2014), 711--729.
    
\bibitem{J1} (MR2257424 (2008f:37177)) [10.3934/dcds.2007.17.159]
    \newblock N. Ju,   
    \newblock The global attractor for the solutions to the 3{D} viscous primitive equations,
    \newblock \emph{Discrete Contin. Dyn. Syst.}, \textbf{17} (2007), 159--179.    
    
\bibitem{JT} (MR3302126) [10.1007/s00332-014-9223-8]
    \newblock N. Ju and R. Temam,   
    \newblock Finite dimensions of the global attractor for 3{D} primitive equations with viscosity,
    \newblock \emph{J. Nonlinear Sci.}, \textbf{25} (2015), 131--155.        
    
\bibitem{daJT} (MR1195597) [10.1016/0167-2789(92)90233-D]
    \newblock Don A. Jones and E. S. Titi,  
    \newblock Determining finite volume elements for the {$2$}{D} {N}avier-{S}tokes equations,
    \newblock \emph{Phys. D}, \textbf{60} (1992), 165--174, Experimental mathematics: computational issues in nonlinear science (Los Alamos, NM, 1991).   

\bibitem{BGMR1} (MR2025211)
    \newblock D. Bresch, F. Guill\'en-Gonz\'alez, N. Masmoudi, and M. A. Rodr\'iguez-Bellido,
    \newblock  Uniqueness of solution for the 2{D} primitive equations with friction condition on the bottom, 
    \newblock  in \emph{Seventh {Z}aragoza-{P}au {C}onference on {A}pplied {M}athematics and {S}tatistics ({S}panish) ({J}aca, 2001)} (eds. E.H. Zarantonello and Author 2), Univ. Zaragoza, \textbf{27} (2003), 135--143.

\bibitem{HTZ} (MR2172554)
    \newblock C. Hu, R. Temam, and M. Ziane,
    \newblock Regularity results for linear elliptic problems related to the primitive equations, 
    \newblock in \emph{Frontiers in mathematical analysis and numerical methods}, World Sci. Publ., (2004), 149--170. 
    
\bibitem{STT} (MR2160019) [10.1007/0-387-24276-7\_60]
    \newblock E. Simonnet, T. Tachim-Medjo, and R. Temam,
    \newblock Higher order approximation equations for the primitive equations of the ocean, 
    \newblock in \emph{Variational analysis and applications}, Springer, \textbf{79} (2005), 1025--1048. 
    
\bibitem{rD}
     \newblock R. Daley,
     \newblock \emph{Atmospheric {D}ata {A}nalysis},
     \newblock Cambridge Atmospheric and Space Science Series, 1993.

\bibitem{eK}
     \newblock E. Kalnay,
     \newblock \emph{Atmospheric {M}odeling, {D}ata {A}ssimilation and {P}redictability},
     \newblock Cambridge University Press, 2003.

\bibitem{LSZ}
     \newblock K. Law, A. Stuart, and K. Zygalakis,
     \newblock \emph{A {M}athematical {I}ntroduction to {D}ata {A}ssimilation},
     \newblock Vol. 62 of Texts in Applied Mathematics, Springer, Cham, 2015.

\bibitem{jP}
     \newblock J. Pedlosky,
     \newblock \emph{Geophysical {F}luid {D}ynamics},
     \newblock Springer New York, 1987.

\bibitem{CF}
     \newblock P. Constantin and C. Foias,
     \newblock \emph{Navier-{S}tokes {E}quations},
     \newblock University of Chicago Press, Chicago, IL, 1988.
     
\bibitem{T}
     \newblock R. Temam,
     \newblock \emph{Navier-{S}tokes {E}quations: {T}heory and {N}umerical {A}nalysis},
     \newblock AMS Chelsea Publishing, Providence, RI, 2001. 

\bibitem{FLT}
\newblock A. Farhat, E. Lunasin, and E. S. Titi,
\newblock On the {C}harney Conjecture of Data Assimilation Employing Temperature
	Measurements Alone: The Paradigm of 3D Planetary Geostrophic Model, 
\newblock {arXiv:1608.04770}.

\bibitem{LP}
\newblock A. Larios and Y. Pei,
\newblock Nonlinear continuous data assimilation, 
\newblock {arXiv:1703.03546}.

\bibitem{MT}
\newblock C. Mondaini, C. and E. S. Titi,
\newblock Postprocessing {G}alerkin method applied to a data assimilation algorithm: a uniform in time error estimate, 
\newblock {arXiv:1612.06998}.

\bibitem{J2}
\newblock N. Ju,
\newblock Global Uniform Boundedness of Solutions to viscous 3D Primitive Equations with Physical Boundary Conditions, 
\newblock {arXiv:1710.04622v2}.

\end{thebibliography}
\end{document}